# Failure of parametric H-principle for maps with prescribed Jacobian

by Joseph Coffey

**Abstract**

Let $M$ and $N$ be closed $n$-dimensional manifolds, and equip $N$ with a volume form $\sigma$. Let $\mu$ be an exact $n$-form on $M$. Arnold then asked the question: When can one find a map $f\colon M \to N$ such that $f^*\sigma = \mu$? In 1973 Eliashberg and Gromov showed that this problem is, in a deep sense, trivial: It satisfies an $h$-principle, and whenever one can find a bundle map $f_{bdl}\colon TM \to TN$ which is degree 0 on the base and such that $f_{bdl}{}^*(\sigma) = \mu$ one can homotop this map to a solution $f$. That is if the naive topological conditions are satisfied on can find a solution. There is no further interesting geometry in the problem.

We show the corresponding parametric $h$-principle fails -- if one considers *families* of maps inducing $\mu$ from $\sigma$, one can find interesting topology in the space $\mathcal{M}_\mu$ of solutions which is not predicted by an $h$-principle. Moreover the homotopy type of such maps is "quantized": for certain families of forms homotopy type remains constant, jumping only at discrete values.

## 1 Introduction

Let $M$ and $N$ be compact $n$-manifolds. Endow $N$ with a volume form $\sigma$, and let $\mu$ be an $n$-form on $M$, We consider $\mathcal{M}_\mu$ - the smooth maps $f\colon M \to N$ whose Jacobian is prescribed by:

$$f^*(\mu) = \sigma$$

In the late sixties Arnold asked when one can find such a map $f$ with preescribed Jacobian. In 1973 Eliashberg and Gromov gave a complete answer to this question; they showed that under some obvious necessary conditions the problem satisfies an $h$-principle. This reduces the geometric question of the existence of an inducing map to the algeraic topology of certain bundle maps, and thus there is no geometry in the *existence* of maps in $\mathcal{M}_\mu$. However in this paper we will show *that there is interesting topology in the space $\mathcal{M}_\mu$, beyond that predicted by an $h$-principle*. Maps with prescribed Jacobian thus lie on the border of soft and hard mathematics, blessed with both many examples and interesting structure.

**What is an "$h$-principle?"** To understand this phrase we consider an easier question than Arnold's: Instead of seeking a map $f\colon M \to N$ inducing $\mu$, we look only for a formal solution, a bundle map

$$(f_{\text{base}}, f_{\text{fb}})\colon TM \to TN$$

such that $f_{fb}^*\sigma = \mu$. Call the space of such formal solutions $\mathcal{B}\mathcal{M}_\mu$. Finding maps in $\mathcal{B}\mathcal{M}_\mu$ is a far more flexible problem, and it is one of standard algebraic topology. There is a natural inclusion of the genuine solutions $\mathcal{M}_\mu$ into these formal solutions $\mathcal{B}\mathcal{M}_\mu$:

$$\begin{aligned} i\colon \mathcal{M}_\mu &\hookrightarrow \mathcal{B}\mathcal{M}_\mu \\ i(g) &\to (g, Dg) \end{aligned}$$

and so finding a formal solution in $\mathcal{B}\,\mathcal{M}_\mu$ is a necessary first step to finding a map with prescribed Jacobian.

Eliashberg and Gromov showed that (under the conditions prescribed below in Theorem 1.1) this necessary step is sufficient; if one has a formal solution $(g,\ g_\text{fib}) \in \mathcal{B}\,\mathcal{M}_\mu$ inducing $\mu$ one can homotop it, through bundle maps in $\mathcal{B}\,\mathcal{M}_\mu$, to a genuine solution $(f,\ D\,f) \in \mathcal{M}_\mu$, where the map along the fibers is actually the derivative of the map on the base. In other words, the induced map on $\pi_0$:

$$i_* \colon \pi_0(\mathcal{M}_\mu) \to \pi_0(\mathcal{B}\,\mathcal{M}_\mu)$$

is surjective. Because one can always homotop a formal "bundle map" solution to a genuine solution, one says that: "Maps with prescribed Jacobian satisfy a homotopy principle or **h-principle.**"

**Theorem 1.1. (h-principle for maps with prescribed Jacobian -- Eliashberg-Gromov-[GÈ73])**
*Let $M$ and $N$ be oriented, compact $n$-manifolds. Let $\sigma$ be a volume form on $N$ and let $\mu$ be an $n$-form on $M$ which alternates in sign. Further suppose that $\int_M \mu = k \cdot \int_N \sigma$ for some $k \in \langle z \rangle$. Then every bundle map $(g,\ g_\text{fib}) \in \mathcal{B}\,\mathcal{M}_\mu$ such that the induced map $g\colon M \to N$ of base spaces is degree $k$, can be homotoped through $\mathcal{B}\,\mathcal{M}_\mu$ to a map $f \in \mathcal{M}_\mu$.*

An $h$-principle is characteristic of soft geometry: immersions, submersions, symplectic immersions symplectic embeddings in codimension $> 2$, Lagrangian immersions, and Legendrian embeddings all satisfy h-principles. All of these problems seem geometric on their face. However they each satisfy an $h$-principle, and thus solving them reduces to the appropriate theory of bundle maps. This is cause for both celebration and mourning, depending on ones point of view. To illustrate consider the case of immersions: If one seeks an immersion of one manifold $M$ into another $N$, the corresponding $h$-principle often reduces the question to one of characteristic classes, thus simplifying the problem considerably. However, if ones goal is to understand the topology of $M$ and $N$, the $h$-principle for immersions says that we will learn nothing new from studying immersions $M \to N$ beyond the relatively coarse theory of the vector bundles $TM$ and $TN$.

Faced with an $h$-principle, the seeker of geometry might then look for structure within the topology of the space of solutions. Here again one has an $h$-principle to contend with: If the inclusion $i$ of formal solutions into genuine ones is a *weak homotopy equivalence* one says that the problem satisfies a **parametric h-principle.** Again, the topology of the space of solutions has interesting geometry precisely when the parametric principle fails.

The vast majority of problems satisfying an $h$-principle also satisfy a parametric $h$-principle. This is no accident; most $h$-principles are (sometimes involved) corollaries of a few major theorems of Gromov [Gro86], Eliashberg and Mishachev [EM02, EM97, EM00, EM98]. When these theorems apply to yield an $h$-principle they usually provide for the stronger parametric version as well. Differential geometry is, for the most part, divided cleanly into soft problems -- which satisfy all $h$-principles and are thus blessed with many examples but posess no structure, and hard problems -- which posess so much structure that they have almost no solutions. There is a philosophy, due to Gromov, where the most interesting problems in differential geometry are those which somehow straddle this divide.

In this paper we will show via example that maps with prescribed Jacobian manage these difficult acrobatics. That while they satisfy an $h$-principle, and thus have many examples, *they fail the parametric h-principle.* The map:

$$i\colon \mathcal{M}_\mu \to \mathcal{B}\,\mathcal{M}_\mu$$

while surjective on $\pi_0$ (for maps of the appropriate degree), need *not* be a homotopy equivalence. There is topology in the space $\mathcal{M}_\mu$ which is not predicted by the bundle maps $\mathcal{B}\mathcal{M}_\mu$, but is rather of geometric origin. Maps with prescribed Jacobian thus have a character which is similar to Legendrian embeddings. Legendrian embedding in a given isotopy class is a relatively simple matter (it satisfies an $h$-principle), however the space of such embeddings is quite interesting (the corresponding inclusion into the space of bundle maps is not a homotopy equivalence). This complicated topology has led to a great deal of interesting mathematics: Thurston-Benniquin invariant, relative contact homology etc.. One might hope that the study of maps with prescribed Jacobian to prove similarily interesting.

## 2 Counterexample to parametric $h$-principle

We will now give an example demonstrating the failure of the parametric $h$-principle for maps with prescribed Jacobian.

**Notation 2.1.** *Let $M$, $N$ denote two copies of $S^2$, and let $\sigma$ be a $C^\infty$ volume form on $N$ such that*

$$\int_N \sigma = 1$$

In what follows $M$ will be the domain of our maps, and $N$ will be the range.

A parametric $h$-principle predicts that the homotopy type of the maps $\mathcal{M}_\mu$, inducing $\mu$ from $\sigma$, should remain the same under certain deformations of the form $\mu$ on the domain. We will show that the homotopy type of $\mathcal{M}_\mu$ is not stable in this (to be prescribed) sense.

To this end we now describe a family of forms $\mu_\kappa$ on the domain $M$.

**Definition 2.2.** *Divide the domain $M = S^2$ into two open hemispheres $H_+$ and $H_-$ along a simple closed curve $\gamma$. Let $\mu_\kappa$ be a family of $C^l$ forms, $l \geq q \in \mathbb{N}$, for $0 < \kappa \leq 1$ on $S^2$ such that:*

1. $\mu_\kappa|_x = 0$ *if and only if $x \in \gamma$. Moreover, if we give coordinates*

   $$(y, \theta) : x \in (-\varepsilon, \varepsilon), \theta \in S^1$$

   *to a neighborhood of $\gamma$ such that $\gamma = \{(y, \theta) : y = 0\}$, $\mu_\kappa = y^{2q+1} dy \wedge d\theta$.*

2. *The total area of each hemisphere satisfies:*

   $$\int_{H_+} \mu_\kappa = +\kappa$$
   $$\int_{H_-} \mu_\kappa = -\kappa$$

When we consider maps of higher smoothness ($C^q$ with $q > 1$) we have to augment the definition of the $h$-principle given in the introduction slightly in order to consider the constraints on the $q$-jet of the map $f$ imposed by the decay of the form $\mu$ near $\gamma$. Denote by $J^q(M, N)$ the space of $q$-jets of $C^q$ maps $M \to N$. Then there is a natural projection $\pi \colon J^q(M, N) \to M$.

**Definition 2.3.** *Let $M = N = S^2$. We consider the following spaces of maps:*

1. $\mathcal{B}\mathcal{M}$: *the space of sections $f \colon M \to J^q(M, N)$ of the projection $\pi$, such that:*

   a. *There is a neighborhood $U_f$ of $\gamma$ such that for $x \in U_f$, each $q$-jet $f(x)$ is the $q$-jet of the germ of a map $g_x$ (which varies with $x$, $f$) such that $g_x^*(\sigma) = \mu_1(x)$.*

   b. *For $x \in H_+$ the 1-jet $f_1 \colon T(M) \to T(N)$ is an orientation preserving linear isomorphism. For $x \in H_-$ the 1-jet $f_1 \colon T(M) \to T(N)$ is an orientation reversing linear isomorphism.*

   c. *The 0-jet $f_0 \colon M \to N$ is a degree 0 map.*

2. $\mathcal{M}$: *the space of degree 0, $C^q$ maps $M \to N$ whose $q$-jets lie in $\mathcal{B}\mathcal{M}$.*

3. *Let $\kappa \in (0, 1]$ then we consider:*

   a. $\mathcal{B}\mathcal{M}_{\mu_\kappa}$: *those sections $f \in \mathcal{B}\mathcal{M}$ such that for every $x \in M$ the $q$-jet $f(x)$ is the $q$-jet of the germ of a map $g_x$ (which varies with $x$, $f$) such that $g_x^*(\sigma) = \mu_\kappa(x)$.*

b. $\mathcal{M}_{\mu_\kappa}$: *the space of degree 0, $C^q$ maps $M \to N$ whose $q$-jets lie in $\mathcal{BM}_{\mu_\kappa}$, i.e. are such that*

$$f^*\sigma = \mu_\kappa$$

*Denote by $j_\kappa : \mathcal{BM}_{\mu_\kappa} \to \mathcal{BM}$ the natural inclusions.*

Note that for each $f \in \mathcal{M}$, $\mathrm{im}(f) = f(H_+) = f(H_-)$.

**Remark 2.4. (Regularity: $l$ and $q$)** In Definitions 2.2 and 2.3 there are two important natural numbers $q$ and $l$. $l$ gives the number of derivatives of the forms $\mu_\kappa$ induced on the domain, and $q$ those of the inducing maps $\mathcal{M}_{\mu_\kappa}$. $l$ is then greater than $q$. Note that $q$, the smoothness of the inducing maps $\mathcal{M}_{\mu_\kappa}$, comes subtly into the definition of the forms $\mu_\kappa$ we aim to induce. This is because the maps we construct in Section Kne26 to provide obstructions to the parametric $h$-principle collapse $\gamma$ and then "rotate one hemisphere against the other". These collapsed and rotated maps are $C^q$ along $\gamma$ only when their $q$-jets vanish there. Thus they induce forms whose $q-1$ jets vanish along $\gamma$ as well.

For this collapsing family to be continuous in $C^q$ it seems that we are further constrained, and we must use maps *whose $2q + 2$ jets vanish along $\gamma$*. This is more delicate and possibly a technical hypothesis. It is possible that if instead the form $\mu$ to be induced is suitably transverse at $\gamma$, then maps inducing $\mu$ do satisfy a parametric $h$-principle. Determining the level required is an interesting question. (See Section Kne26 for further discussion.)

## 2.1 Strategy for finding obstructions to a parametric $h$-principle

**Convention 2.5.** In what follows the value $\kappa = 1$ will play a special role. Thus henceforth, unless otherwise specified we will denote by $\kappa$ a real number confined to the *open interval* $(0, 1)$. When we wish to refer the special parameter value $\kappa = 1$ we will simply write 1.

We have the following diagram of inclusions between the spaces of Definition 2.3:

$$
\begin{array}{ccc}
\mathcal{M}_{\mu_1} & \underset{i_1}{\rightrightarrows} & \mathcal{BM}_{\mu_1} \\
 & & \downarrow{j_1} \\
 & & \mathcal{BM} \qquad 0 < \kappa < 1 \\
 & & \uparrow{j_\kappa} \\
\mathcal{M}_{\mu_\kappa} & \underset{i_\kappa}{\rightrightarrows} & \mathcal{BM}_{\mu_\kappa}
\end{array}
$$

The maps $j_1 : \mathcal{BM}_{\mu_1} \to \mathcal{BM}$, and $j_\kappa : \mathcal{BM}_{\mu_\kappa} \to \mathcal{BM}$ are both homotopy equivalences. For near $\gamma$, $\mathcal{BM}_{\mu_\kappa}$ and $\mathcal{BM}$ each place the same restrictions on $q$-jets. Away from $\gamma$, let $\lambda_x = \frac{\mu_\kappa}{f(x)^*(\sigma)}$. Then

$$f(x) \to (t\lambda_x + 1 - t)f(x)$$

gives a deformation retraction of $j_\kappa$. An analogous retraction, which substitutes $\mu_1$ for $\mu_\kappa$ in the definition of $\lambda_x$, exists for $j_1$.

The parametric $h$-principle predicts that the natural maps $i_1$ and $i_\kappa$, which carry a map to its $q$-jet, are homotopy equivalences. Naively, one might produce obstructions to the $h$-principle by providing a map

$$i_{ob} \colon \mathcal{M}_{\mu_\kappa} \to \mathcal{M}_{\mu_1}$$

such that the diagram:

$$
\begin{array}{ccc}
\mathcal{M}_{\mu_1} & \underset{i_1}{\rightrightarrows} & \mathcal{BM}_{\mu_1} \\
 & & \downarrow{j_1} \\
\uparrow{i_{ob}} & & \mathcal{BM} \qquad 0 < \kappa < 1 \\
 & & \uparrow{j_\kappa} \\
\mathcal{M}_{\mu_\kappa} & \underset{i_\kappa}{\rightrightarrows} & \mathcal{BM}_{\mu_\kappa}
\end{array}
$$

commutes up to homotopy, but the map $i_{ob}$ is not a homotopy equivalence. However it is difficult to produce such a map directly. Instead we will show that:

**Theorem 2.6.** *There are spaces $\mathcal{NS} \subset \widehat{\mathcal{NS}}$ such that:*

1. *For each $0 < \kappa < 1$ there is a homotopy equivalence $\alpha_\kappa \colon \mathcal{M}_{\mu_\kappa} \to \mathcal{NS}$.*

2. *There is a homotopy equivalence $\alpha_1 \colon \mathcal{M}_{\mu_1} \to \widehat{\mathcal{NS}}$.*

3. *The resulting diagram:*

$$
\begin{array}{ccccc}
\widehat{\mathcal{NS}} & \underset{\alpha_1}{\leftarrow} & \mathcal{M}_{\mu_1} & \underset{i_1}{\rightarrow} & \mathcal{B}\,\mathcal{M}_{\mu_1} \\
& & & & \downarrow j_1 \\
\uparrow i_{ob} & & & & \mathcal{B}\,\mathcal{M} \qquad 0 < \kappa < 1 \\
& & & & \uparrow j_\kappa \\
\mathcal{NS} & \underset{\alpha_\kappa}{\leftarrow} & \mathcal{M}_{\mu_\kappa} & \underset{i_\kappa}{\rightarrow} & \mathcal{B}\,\mathcal{M}_{\mu_\kappa}
\end{array}
$$

   *commutes up to homotopy.*

4. *The inclusion $i_{ob} \colon \mathcal{NS} \to \widehat{\mathcal{NS}}$ is not a homotopy equivalence.*

*Thus the parametric h-principle for maps with prescribed Jacobian cannot hold for all forms $\mu_\kappa$ for $0 < \kappa \le 1$*

Note Theorem 2.6 will also show that, at least in this example, the failure of the parametric h-principle is "quantized", for we have the following immediate corollary:

**Corollary 2.7.** *The spaces $\mathcal{M}_{\mu_\kappa}$ are all homotopy equivalent for $0 < \kappa < 1$.*

### 2.1.1 Conventions

Homotopy equivalence in this paper means weak homotopy equivalence. We will often discuss families of maps. For example we may consider a family of diffeomorphisms of $N$:

$$\rho \colon D^k \to \mathrm{Diff}(N)$$

When parameterizing these families, $D^k$ will always refer to the closed $k$-disc, $S^1$ to the circle, and $I$ to the closed interval $[0, 1]$, unless otherwise qualified. When we wish to refer to the map $\rho(d)$ we will write it as $\rho_d$, in order to keep the number parentheses to a minimum.

A **deformation retraction** of a map of pairs $\rho \colon (D^k, \partial D^k) \to (X, Y)$, is a homotopy $\boldsymbol{\rho} \colon (D^k, \partial D^k) \times I \to (X, Y)$ such that $\boldsymbol{\rho}_{d,t} = \rho_d$, and $\boldsymbol{\rho}_{d,1} \in Y$ for every $d \in D^k$, and further $\boldsymbol{\rho}_{d,t} = \rho_d$ for all $d \in \partial D^k$. We say that an inclusion $i \colon X \hookrightarrow Y$ admits a deformation retraction, if every such family $\rho$ does.

# 3 Construction of Model Spaces $\mathcal{NS}$ and $\widehat{\mathcal{NS}}$

We construct the commutative diagram described in the statement of Theorem 2.4.

**Definition 3.1. (Model Spaces)** *We say map $f \in \mathcal{M}$ **overlaps** if $f|_{H_+}$ is not injective. We consider the following subspaces of $\mathcal{M}$:*

1. $\mathcal{NS}_{\mathrm{nc}} : f \in \mathcal{M}$ *which are non-surjective.*

2. $\mathcal{NS} : f \in \mathcal{M}$ *which are non-surjective and overlap.*

3. $\widehat{\mathcal{NS}} : f \in \mathcal{M}$ *such that either $f \in \mathcal{NS}$, or $f$ is surjective but does not overlap.*

*Then denote the inclusion by:*

$$i_{\mathrm{ob}} \colon \mathcal{NS} \hookrightarrow \widehat{\mathcal{NS}}$$

The $\mathcal{NS}$ stands for non surjective, the nc for no (overlapping) constraint. We have the following diagram of inclusions:

$$
\begin{array}{ccc}
\mathcal{M}_{\mu_\kappa} & \underset{i_\kappa}{\to} & \mathcal{NS}_{\mathrm{nc}} \\
& & \uparrow i_{\mathrm{nc}} \\
& & \mathcal{NS} \qquad 0 < \kappa < 1 \\
& & \downarrow i_{\mathrm{ob}} \\
\mathcal{M}_{\mu_1} & \underset{i_1}{\to} & \widehat{\mathcal{NS}}
\end{array}
$$

**Remark 3.2. (On overlapping constraints)**  In subsubsection 3.1.1 we will show that the natural inclusions: $i_\kappa \colon \mathcal{M}_\kappa \hookrightarrow \mathcal{NS}_{\mathrm{nc}}$ and $i_1 \colon \mathcal{M}_1 \hookrightarrow \widehat{\mathcal{NS}}$ are deformation retracts by composition with diffeomorphisms. When $\kappa = 1$ these methods fail if the map does not overlap. Thus we eliminate such maps in our model $\widehat{\mathcal{NS}}$ by fiat. However we pay a price: maps in $\mathcal{M}_\kappa$ need not overlap. Thus they need not lie in $\mathcal{NS}$, but lie only in the larger space $\mathcal{NS}_{\mathrm{nc}}$. We show that this distinction is irrelvant in homotopy: in 3.0.2 we show that we can deform families of maps in $\mathcal{NS}_{\mathrm{nc}}$ to overlap themselves, i.e. that *the inclusion* $i_{\mathrm{nc}} \colon \mathcal{NS} \hookrightarrow \mathcal{NS}_{\mathrm{nc}}$ *is a deformation retract.*

If we denote a corresponding homotopy inverse of $i_{\mathrm{nc}}$ by $r_{\mathrm{nc}}$ the homotopy equivalences required of Theorem 2.6 are given by:

$$
\begin{aligned}
\alpha_\kappa &= r_{\mathrm{nc}} \circ i_\kappa \\
\alpha_1 &= i_1
\end{aligned}
$$

Thus, once we have shown that $i_1$, $i_\kappa$ (subsubsection 3.1.1) and $i_{\mathrm{nc}}$ (subsubsection 3.0.2) are homotopy equivalences we will have established the first portion of Theorem 2.6, namely:

**Proposition 3.3. (Model Spaces are homotopy equivalent)**

1. *For $0 < \kappa < 1$ there is a homotopy equivalence $\alpha_\kappa \colon \mathcal{M}_{\mu_\kappa} \to \mathcal{NS}$*

2. *There is a natural inclusion $\alpha_1 \colon \mathcal{M}_{\mu_1} \to \widehat{\mathcal{NS}}$, and this map is a deformation retract.*

**Proof.**  The proof of this Proposition will occupy the remainder of this section.

### 3.0.2  Making families of maps overlap

In this subsubsection we show:

**Lemma 3.4.**  *The inclusion $i_{\mathrm{nc}} \colon \mathcal{NS} \hookrightarrow \mathcal{NS}_{\mathrm{nc}}$ is a homotopy equivalence.*

**Proof.**  For the proof of this proposition, and that of several other facts in this paper, we will use the following discretization process to reduce to simpler maps:

**Definition 3.5. (Nets of intervals)** *Parametrize $\gamma$ by the unit circle. For each $i \in \mathbb{N}$, let $\mathcal{N}^i$ denote the set of disjoint closed intervals in $\gamma$ which have radius $\frac{2\pi}{10 \cdot 2^i}$ and are centered at $\frac{j 2\pi}{2^i}$ for indexes $1 \le j \le 2^i$ in $\mathbb{N}$. Let:*

1. $\mathcal{NS}^i = \{ f \in \mathcal{NS} \colon$ *There exists a $\Delta_j \in \mathcal{N}^i$ such that $f$ immerses $\Delta_j \}$*

2. $\mathcal{NS}^i_{\mathrm{nc}} = \{ f \in \mathcal{NS}_{\mathrm{nc}} \colon$ *There exists a $\Delta_j \in \mathcal{N}^i$ such that $f$ immerses $\Delta_j \}$*

*For each such interval $\Delta_j \in \mathcal{N}^i$ let:*

1. $\mathcal{U}^j_{\mathrm{nc}} = \{ f \in \mathcal{NS}^i_{\mathrm{nc}} \colon f$ *immerses $\Delta_j \}$*

2. *Let $\mathcal{U}^j = \mathcal{U}^j_{\mathrm{nc}} \cap \mathcal{NS}^i$.*

*These form open coverings of $\mathcal{NS}^i_{\mathrm{nc}}$ and $\mathcal{NS}^i$ respectively.*

The intervals in $\mathcal{N}^i$ satisfy two important properties:

1. **(Intervals are scattered and small)** For each $r > 0$ there is an $i_0 \in \mathbb{N}$ such that for any $i > i_0$, any interval of length $r$ inside $\gamma$ contains some $\Delta \in \mathcal{N}^i$.

2. **(Intervals are nested)** Each interval $\Delta \in \mathcal{N}^i$ contains an interval $\Delta' \in \mathcal{N}^j$ for all $j > i$.

As a result the spaces $\mathcal{NS}^i$ and $\mathcal{NS}^i_{\mathrm{nc}}$ form direct systems under inclusion:

$$\mathcal{NS}^1 \subset \mathcal{NS}^2 \subset ... \mathcal{NS}^i \subset ...$$
$$\mathcal{NS}^1_{\mathrm{nc}} \subset \mathcal{NS}^2_{\mathrm{nc}} \subset ... \mathcal{NS}^i_{\mathrm{nc}} \subset ...$$

Further, for any compact family $\rho : D^k \to \mathcal{NS}$ (or $\mathcal{NS}_{\mathrm{nc}}$) there is an $i_0$ such that $\rho(D^k) \subset \mathcal{NS}^i$ (or $\mathcal{NS}^i{}_{\mathrm{nc}}$) for $i > i_0$.

Thus to prove Proposition 3.4 it is sufficient to show that the inclusion $i_{\mathrm{nc}} : \mathcal{NS}^i \hookrightarrow \mathcal{NS}^i_{\mathrm{nc}}$ is a homotopy equivalence. We now prove this by applying the following elementary (but not so well known) Lemma in homotopy theory to the inclusion $i_{\mathrm{nc}} : \mathcal{NS}^i \hookrightarrow \mathcal{NS}^i_{\mathrm{nc}}$ and the coverings $\mathcal{U}^j, \mathcal{U}^j_{\mathrm{nc}}$:

**Lemma 3.6. (Homotopy Decomposition Lemma [Gra75] - Proposition 16.24)** *Let $f : X \to Y$ be a continuous map. Let $U^j_Y$ be a finite covering of $Y$ by open sets, and denote $f^{-1}(U^j_Y)$ by $U^j_X$. Suppose that for each $J \subset I$ the the restriction*

$$f : \bigcap_{j \in J} U^j_X \to \bigcap_{j \in J} U^j_Y$$

*is a homotopy equivalence then $f$ is a homotopy equivalence*

Note that Proposition 16.24 given by Gray in [Gra75] refers to excisive covers. These are covers by sets whose interiors cover. In particular open covers are excisive. Gray also only covers the case of a covering by 2 sets, but the general case follows by induction. This induction is relatively straightforward, however a sketch is given in the Appendix (section 6.2) for the readers convenience.

It is enough then to show that the map inclusion of each multi intersection:

$$i_{\mathrm{nc}} : \bigcap_{j \in J} \mathcal{U}^j \to \bigcap_{j \in J} \mathcal{U}^j_{\mathrm{nc}}$$

is a homotopy equivalence. We will then apply the Homotopy Decomposition Lemma. These consist of those maps which immerse the union of intervals $\Delta = \bigcup_{j \in J} \Delta_j$.

I claim that each such restriction of $i_{\mathrm{nc}}$ is a deformation retract. Consider a map of pairs

$$\rho : (D^k, \partial D^k) \to (\bigcap_{j \in J} \mathcal{U}^j, \bigcap_{j \in J} \mathcal{U}^j_{\mathrm{nc}})$$

We will now construct a deformation retraction

$$\boldsymbol{\rho} : D^k \times I \to \bigcap_{j \in J} \mathcal{U}^j_{\mathrm{nc}}$$

of $\rho$. We can find a closed interval $\gamma_{\mathrm{em}} \subset \Delta \subset \gamma$ such that $\rho_d|_{\gamma_{\mathrm{em}}}$ is an embedding for all $d \in D^k$. For, $\Delta$ is immersed by every map $\rho_d$, and the derivatives of the maps $\rho_d$ are uniformly bounded as $D^k$ is compact. Then, we apply Lemma 3.7 below to find a homotopy $\boldsymbol{\rho}' : D^k \times I \to \bigcap_{j \in J} \mathcal{U}^j_{\mathrm{nc}}$:

**Lemma 3.7.** *Let $\gamma_{\mathrm{em}} \subset \gamma$ be a closed interval, let $U \subset D^k$ be an open set, and let $\rho : U \to \mathcal{NS}_{\mathrm{nc}}$ is a family of maps which each embedd $\gamma_{\mathrm{em}}$. Then there is a homotopy $\boldsymbol{\rho} : U_\Delta \times I \to \mathcal{NS}_{\mathrm{nc}}$ such that:*

1. $\boldsymbol{\rho}_{d,0} = \rho_d$

2. $\boldsymbol{\rho}_{d,1}$ *overlaps, and thus $\boldsymbol{\rho}_{d,1} \in \mathcal{NS}$.*

3. *If $\boldsymbol{\rho}_{d,0}$ overlaps, then $\boldsymbol{\rho}_{d,t}$ overlaps for all $t \in [0,1]$.*

4. *If $\boldsymbol{\rho}_{d,0}$ immerses an interval $\Delta \subset \gamma$ then $\boldsymbol{\rho}_{d,t}$ immerses $\Delta$ for all $t \in [0,1]$.*

The geometric idea behind the proof of Lemma 3.7 is transparent: one pushes the family of maps $\boldsymbol{\rho}$ along bands whose core's meet the image of $\gamma_{\mathrm{em}}$ transversely as in the figure below. However its proof requires a fair bit of notation and is therefore postponed to the Appendix (see 6.1).

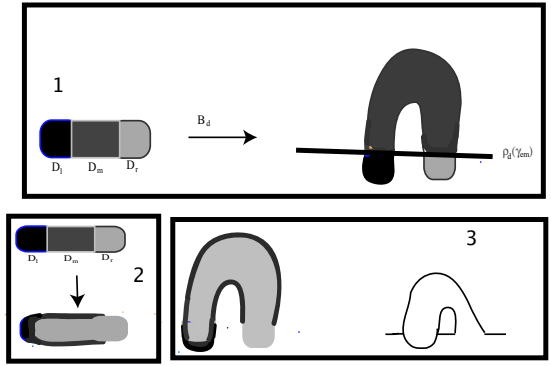

**Figure 3.1.** We prove Lemma 3.7 in the Appendix by first constructing a family of embedded bands $B_d$ which meet $\rho_d(\gamma_{\mathrm{em}})$ in a pair of arcs (step 1 above) Then we push each map along these bands so that they eventually overlap, by conjugating the maps $\rho_d$ with a diffeomorphism of the domain of the embeddings $B_d$ (steps 2 and 3).

To complete the construction of $\boldsymbol{\rho}$ and thus the Proof of Proposition 3.4 we now cut off the deformation $\boldsymbol{\rho}'$ provided by Lemma 3.7 near $\partial D^k$ so that it leaves $\rho$ unchanged there. Note that for each $d \in \partial D^n$, $\rho_d$ overlaps. This is an open condition, and thus we can find a neighborhood $U_\partial$ of $\partial D^k$ where $\rho_d$ overlaps for each $d \in U_\partial$. Then $\boldsymbol{\rho}'_{d,t}$ will also overlap for $d \in U_\partial$, $t \in [0,1]$ (condition 3 of Lemma 3.7). Let $\phi \colon D^k \to [0,1]$ be a smooth function which vanishes on $\partial D^k$ and such that $\phi(d) = 1$ for $d \in D^k \backslash U_\partial$. Then define:

$$\boldsymbol{\rho} \colon D^k \times I \to \bigcap_{j \in J} \mathcal{U}^j_{\mathrm{nc}}$$

by:

$$\boldsymbol{\rho}_{d,t}(x) = \boldsymbol{\rho}'_{d, \phi(d) \cdot t}(x)$$

This completes the proof of Lemma 3.4, modulo that of Lemma 3.7 given in the Appendix (6.1).                    □

## 3.1  The inclusions $i_\kappa \colon \mathcal{M}_{\mu_\kappa} \to \mathcal{NS}_{\mathrm{nc}}$, $i_1 \colon \mathcal{M}_{\mu_1} \to \widehat{\mathcal{NS}}$ are deformation retracts

In the prior subsection we showed that there is a deformation retraction $r_{\mathrm{nc}} \colon \mathcal{NS}_{\mathrm{nc}} \to \mathcal{NS}$. In this subsection we show that the inclusions $i_\kappa \colon \mathcal{M}_{\mu_\kappa} \to \mathcal{NS}_{\mathrm{nc}}$, $i_1 \colon \mathcal{M}_{\mu_1} \to \widehat{\mathcal{NS}}$ admit deformation retracts. Since:

$$\alpha_k = \mathcal{M}_{\mu_\kappa} \underset{i_k}{\to} \mathcal{NS}_{\mathrm{nc}} \underset{r_{\mathrm{nc}}}{\to} \mathcal{NS}$$
$$\alpha_1 = \mathcal{M}_{\mu_\kappa} \underset{i_1}{\to} \mathcal{NS}_{\mathrm{nc}}$$

this will show that $\alpha_k$ and $\alpha_1$ are homotopy equivalences and thus complete the proof of Proposition 3.3.

**Proof.** We begin by showing that the inclusons $i_\kappa \colon \mathcal{M}_{\mu_\kappa} \to \mathcal{NS}_{\mathrm{nc}}$ , $i_1 \colon \mathcal{M}_{\mu_1} \to \widehat{\mathcal{NS}}$ exist and moreover factor through the following intermediate spaces:

**Definition 3.8. (Maps inducing correct area)** *Denote by $\mathcal{V}_1$ the space of volume forms $\omega$ on the range $N$ such that $\int_N \omega = 1$. Let $\omega \in \mathcal{V}_1$. For each $0 < \kappa \le 1$ denote by $(\mathcal{M}_\kappa, \omega)$ the space of maps $f \in \mathcal{M}$ such that $\int_{H_+} f^*\omega = \kappa$. We will abbreviate $(\mathcal{M}_\kappa, \sigma)$ by $\mathcal{M}_\kappa$.*

Then I claim that each inclusion $i_\kappa$, $i_1$ exists and factors as:

$$i_\kappa \colon \mathcal{M}_{\mu_\kappa} \underset{\beta_\kappa}{\overrightarrow{\phantom{xx}}} \mathcal{M}_\kappa \underset{\chi_\kappa}{\overrightarrow{\phantom{xx}}} \mathcal{NS}_{\mathrm{nc}} \quad \text{for } 0 < \kappa < 1$$

$$i_1 \colon \mathcal{M}_{\mu_1} \underset{\beta_1}{\overrightarrow{\phantom{xx}}} \mathcal{M}_1 \underset{\chi_1}{\overrightarrow{\phantom{xx}}} \widehat{\mathcal{NS}}$$

We remind the reader that $\mathcal{NS}_{\mathrm{nc}}$ denotes the nonsurjective maps, and that $\widehat{\mathcal{NS}}$ denotes those which are either non-surjective and overlapping, or surjective and non-overlapping.

Clearly $\mathcal{M}_{\mu_\kappa} \subset \mathcal{M}_\kappa$, and $\mathcal{M}_{\mu_1} \subset \mathcal{M}_1$. To see that there are also natural inclusions: $\chi_\kappa \colon \mathcal{M}_\kappa \hookrightarrow \mathcal{NS}_{\mathrm{nc}}$ and $\chi_1 \colon \mathcal{M}_1 \hookrightarrow \widehat{\mathcal{NS}}$ note that for each $f \in \mathcal{M}$, $f|_{H_+}$ preserves orientation. Thus:

$$\int_{f(H_+)} \sigma \leq \int_{H_+} f^* \sigma$$

with equality if and only if $f|_{H_+}$ is injective. So we see that a map $f \in \mathcal{M}_\kappa$ for $\kappa < 1$ cannot be surjective, and that a map $f \in \mathcal{M}_1$ may be surjective but only when $f|_{H_+}$ is injective, i.e. when it does not overlap.

That $\beta_\kappa$, $\beta_1$ are deformation retracts is an immediate consequence of Moser's Lemma. We now show that the $\chi_\kappa$, $\chi_1$ are also deformation retracts. This will also be via Moser's Lemma, but the details are considerably more involved. They are however similar in spirit to those in [Cof].

### 3.1.1 Strategy to show that $\chi_\kappa$, $\chi_1$ are deformation retracts: Moser's Lemma

For each map $f$ we construct a new volume form $\nu(f)$ on the range $N$ such that

$$f \in (\mathcal{M}_\kappa, \nu(f))$$

and further

$$\int_N \nu(f) = 1$$

Given a disc of maps $\rho$, we thus produce a resulting disc of forms $\nu$. By convexity of volume forms with a given volume, we can contract this disc of forms $\nu$. Moser's Lemma then provides a family of diffeomorphisms inducing this contraction of $\nu$, which, upon postcomposition with the maps $\rho$, moves $\rho$ into $\mathcal{M}_\kappa$ or $\mathcal{M}_1$ respectively. We carry out this program in the remainder of this subsection.

### 3.1.2 Defining a family of forms

**Proposition 3.9.** *Denote by $\mathcal{V}_1$ the volume forms on $N$ of total volume 1. Suppose that either:*

1. *$\rho \colon (D^k, \partial D^k) \to (\mathcal{NS}_{\mathrm{nc}}, \mathcal{M}_\kappa)$ and $\kappa < 1$ or*

2. *$\rho(D^k, \partial D^k) \to (\widehat{\mathcal{NS}}, \mathcal{M}_1)$*

*Then there is a map $\nu \colon (D^k, \partial D^k) \to (\mathcal{V}_1, \sigma)$ $\rho_d \in (\mathcal{M}_\kappa, \nu_d)$ or $\rho_d \in (\mathcal{M}_1, \nu_d)$ respectively.*

**Proof.** We build the family of forms $\nu_d$ by first scaling $\sigma$ in the image of $\rho_d$ so that $\rho_d \in (\mathcal{M}_\kappa, \nu_d)$. Then we scale away from the image so that the resulting forms $\nu_d$ have the proper cohomology class, and thus lie in $\mathcal{V}_1$. To wit:

**Scaling in the image:**

**Lemma 3.10.** *For each $0 < \delta$, there is a continuous function $\phi_{\mathrm{im}}^\delta \colon D^k \times S^2 \to [0, \infty)$ such that each restriction $\phi_{\mathrm{im}}^\delta|_{d \times S^2}$ is $C^\infty$ and $V_d^\delta = (\phi_{\mathrm{im}}^\delta)^{-1}(0) \cap (d \times S^2)$ satisfies:*

1. *$\rho_d^{-1}(V_d^\delta)$ is a closed neighborhood of $\gamma$ such that:*

$$\int_{\rho_d^{-1}(V_d^\delta) \cap H_+} \rho_d^* \sigma < \delta$$

2. *$V_d^\delta \cup \mathrm{im}(\rho_d) = S^2$.*

**Proof.** The proof of this Lemma is straightforward. One firsts constructs a continious functions satisfying conditions 1 and 2. Then one convolutes along the leaves to smooth the resulting function while retaining these conditions. Details are given in subsection 6.3 of the Appendix. □

Now choose $0 < \delta < \kappa$, and apply Lemma 3.10 to obtain a function $\phi^\delta_{\mathrm{im}}$. Now, we apply $\phi^\delta_{\mathrm{im}}$ to scale $\sigma$ in the image of $\rho_d$ so that $\rho_d \in (\mathcal{M}_\kappa, \nu_d)$. Note that the $(1 + \phi^\delta_{\mathrm{im}}(d, x))^s$ converge to the characteristic functions of the sets $V^\delta_d$ as $s \to -\infty$. Thus

$$\lim_{s \to -\infty} \int_{H_+} \rho^*_d ((1 + \phi^\delta_{\mathrm{im}}(d, x))^s \cdot \sigma(x)) = \int_{\rho^{-1}_d(V^\delta_d) \,\cap\, H_+} \rho^*_d(\sigma(x)) < \delta$$

while

$$\lim_{s \to \infty} \int_{H_+} \rho^*_d ((1 + \phi^\delta_{\mathrm{im}}(d, x))^s \cdot \sigma(x)) = \infty$$

is a monotone increasing function of the real parameter $s$, tak of all values in the interval $(\delta, \infty)$. Thus, as $\delta < \kappa$, there is by the intermediate value theorem there some parameter $s$ for which:

$$\int_{H_+} \rho^*_d ((1 + \phi^\delta_{\mathrm{im}}(d, x))^s \cdot \sigma(x)) = \kappa$$

Further, since the above integral is monotone in $s$ and continious in $d$, we can find a continuous function $s \colon D^k \to \mathbb{R}$ such that for

$$\nu'_d(x) = (1 + \phi^\delta_{\mathrm{im}}(d, x))^{s_d} \cdot \sigma(x)$$

$\int_{H_+} \rho^*_d(\nu'_d) = \kappa$ and thus $\rho_d \in (\mathcal{M}_\kappa, \nu'_d)$.

**Scaling away from the image:** We now alter each form $\nu'_d$ within $N \backslash \mathrm{im}(\rho_d)$ so that the resulting forms reside in $\mathcal{V}_1$, the volume forms of total integral one. As before we construct a scaling function $\phi^\varepsilon_{\mathrm{ms}}$, but this time we scale in the complement of $\mathrm{im}(\rho_d)$.

**Lemma 3.11.** *Let* $\varepsilon \colon D^k \to [0, \infty)$ *be a continious function such that* $\varepsilon_d = 0$ *if and only if* $\rho_d$ *is surjective. Then there is a continious function* $\phi^\varepsilon_{\mathrm{ms}} \colon D^k \times S^2 \to [0, \infty)$ *such that:*

1. *Each restriction* $\phi^\varepsilon_{\mathrm{ms}}|_{d \times S^2}$ *is* $C^\infty$.

2. $V^\varepsilon_d = \phi^\varepsilon_{\mathrm{ms}}{}^{-1}(0) \cap (d \times S^2)$ *is a closed neighborhood of* $\mathrm{im}(\rho_d)$ *within* $s \times S^2$ *satisfying:*

$$\int_{V^\varepsilon_d} \nu'_d - \int_{\mathrm{im}(\rho_d)} \nu'_d \le \varepsilon_d$$

**Proof.** The construction of $\phi^\varepsilon_{\mathrm{ms}}$ entirely analogous to that $\phi^\delta_{\mathrm{im}}$ in Lemma 3.10. We provide it in subsection 6.3 of the Appendix. □

We now wish to find find a continuous function $s' \colon D^k \to \mathbb{R}$ such that

$$\int_N (1 + \phi^\varepsilon_{\mathrm{ms}}(d, x))^{s'_d} \cdot \nu'_d(x)) \;=\; 1$$

Then

$$\nu_d(x) = (1 + \phi^\varepsilon_{\mathrm{ms}}(d, x))^{s'_d} \cdot \nu'_d(x)$$

will give a disc of forms $\nu \colon (D^k, \partial D^k) \to (\mathcal{V}_1, \sigma)$. Further, since $\phi^\varepsilon_{\mathrm{ms}}{}^{-1}(0) \cap (d \times S^2)$ is a closed neighborhood of $\mathrm{im}(\rho_d)$ our forms remain unchanged there and so we will still have $\rho_d \in (\mathcal{M}_\kappa, \nu_d)$.

We divide our proof into two cases:

**Case 1:** $\kappa < 1$, $\rho \colon (D^k, \partial D^k) \to (\mathcal{NS}_{\mathbf{nc}}, \mathcal{M}_\kappa)$: Choose

$$\varepsilon \colon D^k \to [0, 1 - \kappa) \tag{3.1}$$

Then Lemma 3.11 provides a function $\phi_{\mathrm{ms}}^{\varepsilon}\colon D^k \times S^2 \to [0,\infty)$ which vanishes on $V_d^{\varepsilon}$ and is positive elsewhere. The integral

$$\int_N \left(1 + \phi_{\mathrm{ms}}^{\varepsilon}(d,x)\right)^{s'} \cdot \nu_d'(x)$$

is thus a monotone function of the parameter $s'$, taking all values in $(\int_{V_d^{\varepsilon}} \nu_d',\, \infty)$. As before we aim to apply the intermediate value theorem. I claim that with this choice of $\varepsilon$, $\int_{V_d^{\varepsilon}} \nu_d' < 1$. For note that:

$$
\begin{aligned}
\int_{V_d^{\varepsilon}} \nu_d' \;&<\; \int_{H_+} \rho_d^*(\nu_d') + \left(\int_{V_d^{\varepsilon}} \nu_d' - \int_{\mathrm{im}(\rho_d)} \nu_d'\right) \\
&<\; \int_{H_+} \rho_d^*(\nu_d') + \varepsilon_d \\
&=\; \kappa + \varepsilon_d \\
&<\; 1
\end{aligned}
$$

by our choice of $\varepsilon_d$ (equation 3.1).

Thus we can again apply the intermediate value theorem to find a continuous function $s'\colon U \to \mathbb{R}$ such that

$$\int \left(1 + \phi_{\mathrm{ms}}^{\varepsilon}(d,x)\right)^{s_d'} \cdot \nu_d'(x) = 1$$

Thus

$$\nu_d(x) = \left(1 + \phi_{\mathrm{ms}}^{\varepsilon}(d,x)\right)^{s_d'} \cdot \nu_d'(x)$$

gives a disc of forms $\nu\colon (D^k, \partial D^k) \to (\mathcal{V}_1, \sigma)$.

**Case 2: $\kappa = 1$, $\rho(D^k, \partial D^k) \to (\widehat{\mathcal{NS}}, \mathcal{M}_1)$**   Again, we seek a continuous function $s'\colon D^k \to \mathbb{R}$ such that

$$\int_N \left(1 + \phi_{\mathrm{ms}}^{\varepsilon}(d,x)\right)^{s_d'} \cdot \nu_d'(x) \;=\; 1$$

Let $U \subset D^k$ denote those parameters $d$ such that $\rho_d$ is not surjective. We will begin by defining $s'|_U$. Since $\rho_d \in \widehat{\mathcal{NS}}$, such a surjective map must overlap. We can quantify the overlapping of a map $f$ by:

**Definition 3.12.** *Let $f \in \mathcal{M}$, and let $\omega$ be a volume form on the range $N$. We define the **overlap $O(f, \omega)$** of a map $f \in \mathcal{M}$ by:*

$$O(f, \omega) = \int_{H_+} f^*\omega - \int_{\mathrm{im}(f)} \omega$$

*Note that $O(f, \omega) \geq 0$, and $O(f, \omega) > 0$ if an only if $f$ overlaps. We abbreviate $O(f, \sigma)$ by $O(f)$.*

Let $\varepsilon\colon D^k \to [0,\infty)$ be given by $\varepsilon_d = \frac{O(\rho_d, \nu_d')}{2}$. Then for $d \in U$:

$$
\begin{aligned}
\int_{V_d^{\mathrm{im}(\varepsilon)}} \nu_d' \;&=\; \int_{H_+} \rho_d^*(\nu_d') - \left(\int_{H_+} \rho_d^* \nu_d' - \int_{\mathrm{im}(\rho_d)} \nu_d'\right) + \left(\int_{V_\varepsilon^d} \nu_d' - \int_{\mathrm{im}(\rho_d)} \nu_d'\right) \\
&=\; \int_{H_+} \rho_d^*(\nu_d') - O(\rho_d, \nu_d') + \left(\int_{V_\varepsilon^d} \nu_d' - \int_{\mathrm{im}(\rho_d)} \nu_d'\right) \\
&=\; 1 - O(\rho_d, \nu_d') + \varepsilon_d \\
&=\; 1 - \frac{O(\rho_d, \nu_d')}{2} \\
&<\; 1
\end{aligned}
$$

The last inequality follows since $\rho_d \in \widehat{\mathcal{NS}}$ is nonsurjective and thus must overlap. Since $\int_{V_d^{\mathrm{im}(\varepsilon)}} \nu_d' < 1$ we can again construct $s'|_U$ by the intermediate value theorem.

Away from $U$, $\rho_d$ is surjective, and $\phi_{\mathrm{ms}}^{\varepsilon}(d,x) = 0$ for all $x$. Thus

$$\int_N \left(1 + \phi_{\mathrm{ms}}^{\varepsilon}(d,x)\right)^{s_d'} \cdot \nu_d'(x) = 1$$

regardless of the parameter $s'$, and so we can extend the function $s'|_U$ in anyway we like to the rest of the disc $D^k$. Thus we again construct the requisite disc of forms $\nu$ by

$$\nu_d = (1 + \phi^\varepsilon_{\mathrm{ms}}(d)^{s'_d}) \cdot \nu'_d$$

This complete the proof of Proposition 3.9.                                                                    □

### 3.1.3 Applying Moser to gain retraction

We remind the reader that $\mathcal{V}_1$ denotes the space of volume forms on $N$ with volume 1. Denote the diffeomorphisms of $N$ by $\mathrm{Diff}(N)$.

As $\mathcal{V}_1$ is convex, one can construct a retraction $\boldsymbol{\nu}\colon D^k \times I \to \mathcal{V}_1$ of $\nu$ to the constant map by:

$$\boldsymbol{\nu}_{d,t} = t\sigma + (1-t)\nu_d$$

Moser's Lemma then provides a family of diffeomorphisms:

$$\varsigma\colon D^k \times [0,1] \to \mathrm{Diff}(N)$$

such that:

1. $\varsigma^*_{d,t}\sigma = \boldsymbol{\nu}_{d,t}$

2. $\varsigma_{d,0} = id$

3. $\varsigma_{d,t} = id$ for $d \in \partial D^k$

Post composition with $\varsigma$ provides the required retraction of $\rho$ into $\mathcal{M}_\kappa$:

$$\boldsymbol{\rho}_{d,t} = \varsigma_{d,t} \circ \rho_d$$                                                    □

This completes the proof of Proposition 3.3, for we have shown that the inclusions:

$$i_\kappa\colon \mathcal{M}_{\mu_\kappa} \underset{\beta_\kappa}{\overrightarrow{\phantom{xx}}} \mathcal{M}_\kappa \underset{\chi_\kappa}{\overrightarrow{\phantom{xx}}} \mathcal{NS}_{\mathrm{nc}} \quad \text{for } 0 < \kappa < 1$$

$$i_1\colon \mathcal{M}_{\mu_1} \underset{\beta_1}{\overrightarrow{\phantom{xx}}} \mathcal{M}_1 \underset{\chi_1}{\overrightarrow{\phantom{xx}}} \widehat{\mathcal{NS}}$$

are both homotopy equivalences. The homotopy equivalences $\alpha_\kappa$, $\alpha_1$ required of Proposition 3.3 are then given by:

$$\alpha_\kappa = r_{\mathrm{nc}} \circ i_\kappa$$
$$\alpha_1 = i_1$$

where $r_{\mathrm{nc}}$ is the homotopy inverse of $i_{\mathrm{nc}}$.                                          □

## 3.2 The diagram of Theorem 2.6 commutes up to homotopy

We now show that diagram of Theorem 2.6 commutes up to homotopy. Denote the deformation retracts of $i_1$ and $i_\kappa$ by $r_1\colon \widehat{\mathcal{NS}} \to \mathcal{M}_{\mu_1}$ and $r_\kappa\colon \mathcal{NS} \to \mathcal{M}_{\mu_\kappa}$. Then it is enough to show that the diagram:

$$
\begin{array}{ccccc}
\widehat{\mathcal{NS}} & \underset{r_1}{\overrightarrow{\phantom{xx}}} & \mathcal{M}_{\mu_1} & \underset{i_1}{\overrightarrow{\phantom{xx}}} & \mathcal{B}\,\mathcal{M}_1 \\
 & & & & \downarrow{j_1} \\
\uparrow{i_{ob}} & & & & \mathcal{B}\,\mathcal{M} \quad 0 < \kappa < 1 \\
 & & & & \uparrow{j_\kappa} \\
\mathcal{NS} & \underset{r_\kappa \circ i_{\mathrm{nc}}}{\overrightarrow{\phantom{xx}}} & \mathcal{M}_{\mu_\kappa} & \underset{i_\kappa}{\overrightarrow{\phantom{xx}}} & \mathcal{B}\,\mathcal{M}_\kappa
\end{array}
$$

commutes up to homotopy. However this is clear; one can add maps $d_1\colon \widehat{\mathcal{NS}} \to \mathcal{B}\,\mathcal{M}$ and $d_\kappa\colon \mathcal{NS} \to \mathcal{B}\,\mathcal{M}$ which carry each map to its $q$-jet. The commutativity of each triangle within the resulting diagram (up to homotopy) then follows from the naturality of the derivative.

# 4   Obstruction cycle

## 4.1   Strategy to show that $i_{\mathrm{ob}}$ is not a homotopy equivalence

Thus far we have shown (or assumed supposing a parametric $h$-principle) that every  map in the diagram:

$$\begin{array}{ccccc}
\widehat{\mathcal{NS}} & \xleftarrow{\ \alpha_1\ } & \mathcal{M}_{\mu_1} & \xrightarrow{\ i_1\ } & \mathcal{BM}_{\mu_1} \\
 & & & & \downarrow j_1 \\
\uparrow i_{ob} & & & & \mathcal{BM} \qquad 0 < \kappa < 1 \\
 & & & & \uparrow j_\kappa \\
\mathcal{NS} & \xleftarrow{\ \alpha_\kappa\ } & \mathcal{M}_{\mu_\kappa} & \xrightarrow{\ i_\kappa\ } & \mathcal{BM}_{\mu_\kappa}
\end{array}$$

is a homotopy equivalence, save the inclusion $i_{ob} \colon \mathcal{NS} \hookrightarrow \widehat{\mathcal{NS}}$, and that the diagram commutes up to homotopy.

We will now show that $i_{ob}$ is not a homotopy equivalence. In subsection 4.2 we describe an element $[h] \in \pi_1(\mathcal{NS}, \widehat{\mathcal{NS}})$. In the remaining subsections we show that $[h]$ is nontrivial. Our loop $h$ will collapse $\gamma$, and then "rotate $H_+$ against $H_-$". We will show that $[h]$ is nontrivial in $\pi_1(\mathcal{NS}, \widehat{\mathcal{NS}})$ by showing:

1. That it induces a loop $\psi$ of homeomorphisms of the disc, which has a nontrivial free homotopy class.

2. That any deformation of $h$ into $\mathcal{NS}$ would induce a contraction of the loop $\psi$.

Together these two statements give a contradiction. Thus their proof will show that $h$ cannot be deformed into $\mathcal{NS}$. Therefore $[h]$ must be nontrivial, and so $i_{ob}$ is not a homotopy equivalence. Statement 1 is proven in this section. Statement 2 is proven in section Kne26 below. Together they complete the proof of Theorem 2.6 and this paper.

## 4.2   Construction of  obstruction cycle  $[h] \in \pi_1(\mathcal{NS}, \widehat{\mathcal{NS}})$

In the remainder of this section we will construct our loop of maps $h \colon S^1 \to \widehat{\mathcal{NS}}$, and show how it induces the loop $\psi$ of homeomorphisms of the disc. We will proceed as follows:

In 4.2.1 we will construct a path of maps $f \colon I \to \mathcal{M}$ which collapses $\gamma$ to a point. Then in 4.2.2 we will use $f$, along with a rotation of the image of $H_+$ against that of $H_-$, to construct a loop

$$h' \colon S^1 \to \mathcal{M}$$

This will morally be our obstruction cycle, and in 4.3 we will show that it induces a loop of homeomorphisms $\psi'$ of the disc. Further, we will show that any homotopy of $h'$ within $\mathcal{M}$ induces a homotopy of the loop $\psi'$.

However the loop $h'$ we construct will not lie within $\widehat{\mathcal{NS}}$ because it will contain maps which are at once both non-surjective and non-overlapping. We remedy this in 4.4. Here we will apply Lemma 3.7 to overlap the loop $h'$ and deform it to a loop $h \colon S^1 \to \widehat{\mathcal{NS}}$, which gives our  $[h] \in \pi_1(\mathcal{NS}, \widehat{\mathcal{NS}})$.  This homotopy of $h'$ to $h$ in turn induces a homotopy of $\psi'$ to a loop of homeomorphisms $\psi$.

### 4.2.1   A collapsing path of maps $f \colon I \to \mathcal{M}$

Let $(r, \theta)$ denote spherical coordinates on $M = N = S^2$. More precisely

1. Let $r \colon S^2 \to [0, 2]$ denote a smooth height function on $S^2$, chosen so that: $r$ has a unique maximum $x_+ \in H_+$ with maximum value $r(x_+) = 2$ and a unique minimum $x_- \in H_-$ with minimum value $r(x_-) = 0$. Further we arrange that $r$ has no other critical points, and $r^{-1}(1) = \gamma$.

2. $\theta \colon S^2 \backslash \{x_+, x_-\} \to S^1$ is a submersion, such that $\theta|_{r=t}$ is a diffeomorphism for $t \in (0, 2)$.

Further, away from a neighbrhood of $x_+$, the coordinates $(r, \theta)$ identify $\sigma$ with the standard area form on a disc in $\mathbb{R}^2$.

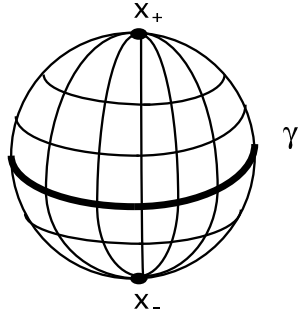

**Figure 4.1.** Level sets of $r$ and $\theta$ respectively.

**Proposition 4.1. (Collapsing path of maps)** *There is a path of maps $f\colon I \to \mathcal{M}$ such that:*

1. *For $t \in [0, 1)$ $f_t$ maps $\gamma$ to $r^{-1}(1-t)$.*

2. *The restriction of $f|_{H_\pm}$ gives a diffeomorphism $f\colon H_\pm \to r^{-1}((1-t, 2])$*

*In particular $f_1$ is such that $f_1(\gamma) = x_-$, and its restriction to either $H_+$ or $H_-$ gives a diffeomorphism $f_1\colon H_\pm \to S^2 \backslash x_-$.*

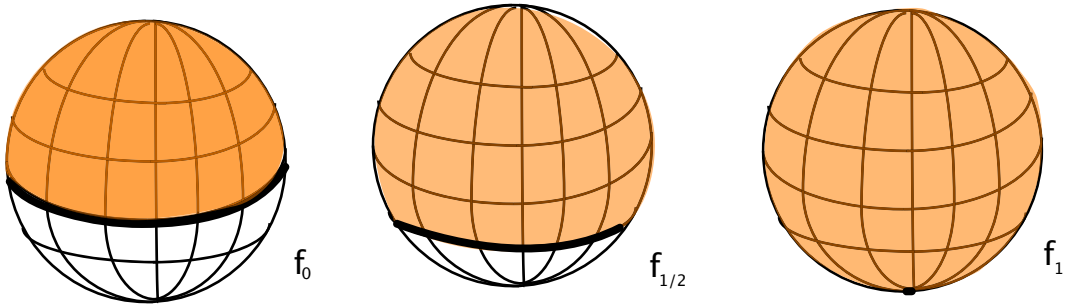

**Figure 4.2.** The collapsing path of maps $f_t$: first $f_0$ folds the sphere along the equator $\gamma$. Then the path $f_t$ collapses $\gamma$ onto the south pole $x_- = r^{-1}(0)$.

**Proof.** Let $f_0\colon S^2 \to S^2$ be given by

$$f_0(r, \theta) = ((1-r)^{2q+2}, \theta)$$

We now show how to extend $f_0$ to the required family of maps $f_t$. The most delicate point is to ensure that as the maps $f_t$ collapse $\gamma$ to $x_-$, we still have that $f_t^*(\sigma) = \mu_1$ near $\gamma$.]

Let $\eta\colon [1, 2] \to [0, 1]$ be a smooth, monotone decreasing function such that $\eta|_{[1, 1\frac{1}{4}]} = 1$, and $\eta|_{[1\frac{3}{4}, 2]} = 0$. Now we define a one parameter family of maps $\varsigma_t\colon \bar{H}_+ \to S^2$ by:

$$\varsigma_t(\theta, r) = (\theta, \sqrt{r^2 - \eta(r)t^2})$$

Then:

1. The map $\varsigma_0$ is the identity.

2. Each $\varsigma_t$ are embeddings for $t < 1$, and $\varsigma_1$ is an embedding away from $\gamma$.

3. Near $\gamma$, $\varsigma_t$ has the form:

$$(\theta, r) \to (\theta, \sqrt{r^2 - t^2})$$

4. $\varsigma_t^*\sigma = \sigma$ near $\gamma$. For away from $x_+$, $\sigma$ is a multiple of the standard area form on a disc in $\mathbb{R}^2$:

$$\sigma = \lambda r\, dr \wedge d\theta \text{ for some } \lambda \in \mathbb{R}$$

$$\begin{aligned}
\varsigma_t^*(\sigma) &= \varsigma_t^*(\lambda r\, dr \wedge d\theta) \\
&= \lambda(\sqrt{r^2 - t^2}\frac{r}{\sqrt{r^2 - t^2}}\, dr \wedge d\theta) \\
&= \lambda r\, dr \wedge d\theta \\
&= \sigma
\end{aligned}$$

Therefore $f_t$, defined by

$$f_t = \varsigma_t \circ f_0$$

gives the requisite collapsing family. For note that as each embedding $\varsigma_t$ preserves $\sigma$ near $\gamma$

$$\begin{aligned}
f_t^*(\sigma) &= f_0^*(\varsigma_t^*(\sigma)) \\
&= f_0^*(\sigma) \\
&= \mu_1
\end{aligned}$$

near $\gamma$. Further, as $\varsigma_1$ collapses $\gamma$ to $x_- = r^{-1}(0)$, so does $f_1 = f_0 \circ \varsigma_1$. Finally, note that $f_t = f_0 \circ \varsigma_t$ is a $C^q$ continuous family of maps, and each map $f_t$ has a vanishing $q$-jet. This is obvious for $t < 1$, at $t = 1$ one computes directly that the map

$$f_1 = \varsigma_1(\theta, (1 - r)^{2q+2})$$

has vanishing derivatives up to order $q$. $\qquad\qquad\qquad\qquad\qquad\qquad\qquad\qquad\qquad\qquad\qquad\square$

### 4.2.2 Collapsing and twisting; our "moral" obstruction cycle $h'$

We now use the collapsing path $f : I \to \mathcal{M}$ to define an element $[h] \in \pi_1(\mathcal{NS}, \widehat{\mathcal{NS}})$. Let $R_\theta$ denote the rotation of $S^2$ about $y$ and its antipode through the angle $\theta$. Initially we will define a loop

$$h' : S^1 \to \mathcal{M}$$

This loop $h'$ will not lie within $\widehat{\mathcal{NS}}$ because it will contain maps which are at once both non-surjective and non-overlapping. Then we will apply Lemma 3.7 to overlap the loop $h'$ and deform it into $\widehat{\mathcal{NS}}$.

**Definition 4.2.** *Denote by $h' : S^1 \to \mathcal{M}$ the following loop of maps:*

1. *$h'|_{[0, \frac{1}{3}]}$ is the "collapsing path" given by $h_t = f_{3t}$*

2. *$h'|_{[1/3, 2/3]}$ is given by*

$$h_t'(x) = \begin{cases} R_{6\pi(t-1/3)} \circ f_1(x) & \text{if } x \in H_+ \\ f_1(x) & \text{if } x \in H_- \end{cases}$$

3. *$h'|_{[2/3, 1]}$ is the "uncollapsing path" given by $h_t' = f_{3-3t}$*

**Remark 4.3.** The loop $h'$ travels through $C^q$ maps because $f_1$ has vanishing $q$-jet along $\gamma$. Note that for the map $f_1$ to have its $q$ jet vanish along $\gamma$ our method requires that the maps $f_t$ have vanishing $2q + 2$ jets. It is an interesting question whether this hypothesis is necessary or an artifact of the method. (see Defintion 2.2 and Remark 2.4).

## 4.3 Obstruction cycle $h'$ induces loop of homeomorphisms $\psi'$

**Definition 4.4.** *Let $\mathcal{H}(H_-, H_+)$ denote the homeomorphisms from $H_-$ to $H_+$, equipped with the topology of convergence on compact sets. Let $\mathcal{H}$ denote the homeomorphisms of the open disc $H_- = D^2$, equipped with same topology.*

**Definition 4.5.** *For each $t \in (0,1)$ and $x \in H_-$ there is a unique $y \in H_+$ such that $h_t'(x) = h_t'(y)$. Denote by $\zeta \colon S^1 \to \mathcal{H}(H_-, H_+)$ the loop of homeomorphisms from $H_-$ to $H_+$ given by:*

$$\zeta_t(x) = y$$

*Define a loop $\psi' \colon S^1 \to O(2) \subset \mathcal{H}$ of homeomorphisms of $H_-$ by $\psi_t' = \zeta_0^{-1} \circ \zeta_t$.*

By considering $\psi' \colon S^1 \to \mathcal{H}$ as loop in $\mathcal{H}$, rather than as a loop of diffeomorphisms with a stronger topology, we discard both control on its derivatives and to some extent its behavior at infinity. Both of these relaxations will be convenient in the constructions of section 5, where assuming that $h$ is homotoped into $\mathcal{NS}$ we contract the loop $\psi$. Finally, note that $[\psi']$ gives a full rotation of the disc, and thus $[\psi'] \in \pi_1(O(2))$ is non-trivial.

## 4.4   Homotoping $h'$ to $h \colon S^1 \to \widehat{\mathcal{NS}}$

We now "overlap" the loop $h'$ so that it lies in $\widehat{\mathcal{NS}}$.

**Proposition 4.6.** *There is a loop $h \colon S^1 \to \widehat{\mathcal{NS}}$ such that:*

1. *$h$ is homotopic to $h'$ within $\mathcal{M}$, i.e. there is a homotopy $\boldsymbol{h} \colon S^1 \times I \to \mathcal{M}$ such that $\boldsymbol{h}_{t,0} = h_t'$ and $\boldsymbol{h}_{t,1} = h_t$.*

2. *$h_0 \in \mathcal{NS}$.*

*The loop $h \colon S^1 \to \widehat{\mathcal{NS}}$ gives an element $[h] \in \pi_1(\mathcal{NS}, \widehat{\mathcal{NS}})$. Moreover, the homotopy $\boldsymbol{h}$ induces a homotopy $\boldsymbol{\zeta} \colon S^1 \times I \to \mathcal{H}(H_-, H_+)$ of $\zeta'$ to a loop $\zeta$. This in turn induces a homotopy $\boldsymbol{\psi} \colon S^1 \times I \to \mathcal{H}$ of $\psi'$ to a loop $\psi$ such that $\psi_t = \zeta_0^{-1} \circ \zeta_t$.*

**Proof.** Let $\gamma_{\mathrm{em}} \subset \gamma$ denote an interval. We now apply Lemma 3.7 to $h, \gamma_{\mathrm{em}}$, and the open set $U \subset S^1 = [0, 1]/\sim$ given by

$$U = S^1 \backslash [\frac{1}{3}, \frac{2}{3}]$$

As each map in $h_U$ embeds $\gamma_{\mathrm{em}}$, Lemma 3.7 provides us a homotopy $\boldsymbol{h}' \colon U \times I \to \mathcal{M}$, to a map $\boldsymbol{h}_{\cdot,1}' \colon U \to \mathcal{NS}$. We wish to extend $\boldsymbol{h}'$ to $\boldsymbol{h} \colon S^1 \times I \to \mathcal{M}$, such that $h = \boldsymbol{h}_{\cdot,1}$ maps $S^1 \to \widehat{\mathcal{NS}}$. However, since surjective maps in $\widehat{\mathcal{NS}}$ must not overlap, we must cut off the overlapping homotopy $\boldsymbol{h}'$ near the boundary of $U$. To do this we introduce quantified notions of surjectivity and overlapping:

We remind the reader that the overlap $O(f)$ of a map $f \in \mathcal{M}$ is defined by:

$$O(f) = \int_{H_+} f^* \sigma - \int_{\mathrm{im}(f)} \sigma$$

We define the **missed area** $MA(f)$ of a map $f \in \mathcal{M}$ by:

$$MA(f) = \int_{S^2 \backslash \mathrm{im}(f)} \sigma$$

Note that a map $f \in \mathcal{M}$ is nonsurjective if and only if $MA(f) > 0$, and overlaps if and only if $O(f) > 0$.

Let $s \colon S^1 \to [0, 1]$ be defined by:

$$s_t \;=\; \inf(\{\tau \colon O(\boldsymbol{h}_{t,\tau}) = MA(h_t')\}, 1)$$

Note that the function $s$ vanishes outside of $U \subset S^1$, as for $t \in S^1 \backslash U$ the maps $h_t$ are surjective. Finally, define $\boldsymbol{h} \colon S^1 \times I \to \mathcal{M}$ by:

$$\boldsymbol{h}_{t,\tau} = \begin{cases} \boldsymbol{h}_{t,s_t \cdot \tau}' & \text{if } t \in U \\ h_t' & \text{otherwise} \end{cases}$$

and define $h \colon S^1 \to \widehat{\mathcal{NS}}$ by:

$$h_t = \boldsymbol{h}_{t,1} \qquad\qquad\qquad\qquad \square$$

**Proposition 4.7.** *The free homotopy class of $\psi\colon S^1 \to \mathcal{H}$ is nontrivial.*

**Proof.** $\psi'$ gives a full rotation of the disc, and thus $[\psi']$ gives a nontrivial element in $\pi_1(O(2))$. The inclusion $O(2) \hookrightarrow \mathcal{H}$ is a deformation retract [Kne26] , and thus $[\psi']$ is also nontrivial in $\pi_1(\mathcal{H})$. Therefore the conjugacy class of $[\psi']$ is nontrivial in $\pi_1(\mathcal{H})$, and so the free homotopy class of $\psi'$ is also nontrivial. Since $\psi$ is freely homotopic to $\psi'$ it lies in this nontrivial class as well. $\qquad\square$

# 5   Homotopy of $h$ into $\mathcal{NS}$ would imply that $\psi$ is contractible

This last section is devoted to the proof of:

**Proposition 5.1.** *If $[h] \in \pi_1(\widehat{\mathcal{NS}}, \mathcal{NS})$ is trivial then the loop $\psi\colon S^1 \to \mathcal{H}$ is freely contractible.*

Note that if $[h] \in \pi_1(\widehat{\mathcal{NS}}, \mathcal{NS})$ is trivial, Propositions 5.1 and 4.7 combine to give a contradiction. Thus the proof of Proposition 5.1 will show that in fact $[h]$ is nontrivial. It will complete the proof of Theorem 2.6, and show the failure of parametric $h$-principle for maps with prescribed Jacobian. It will thus complete this paper.

The proof of Proposition 5.1 is based on the following observation:

> If $h$ admits a homotopy to a loop $\tilde{h}\colon \mathcal{S}^1 \to \mathcal{NS}$ then the maps $\tilde{h}_t$ each miss some open set in the range. Thus, as $\tilde{h}_t(\gamma)$ is the boundary of $\mathrm{im}(\tilde{h}_t)$, each map in $\tilde{h}_t$ must immerse some open interval in $\gamma$.

In the next subsection we will show (after gaining transversality) that the corresponding loop of homeomorphisms $\tilde{\psi}_t$ must fix these immersed intervals. Then in section 5 we show that these fixed intervals allow us to contract the loop $\tilde{\psi}$.

## 5.1   $h$ homotoped into $\mathcal{NS}$ implies $\psi$ can be homotoped to loop $\tilde{\psi}$ where each map fixes some interval $\Delta_j \in \mathcal{N}^i$

**Convention 5.2.** We say that a homeomorphism $f$ of the open disc $D^2$ **fixes** a point $x$ on the boundary if it extends continiously to be the identity map at $x$. We say that $f$ fixes a subset of the boundary if it fixes it pointwise.

This subsection is devoted to the proof of the following Proposition:

**Proposition 5.3.** *Suppose $[h] \in \pi_1(\widehat{\mathcal{NS}}, \mathcal{NS})$ is trivial. Then for $i \gg 0$, $\psi$ can be homotoped within $\mathcal{H}$ to a loop $\tilde{\psi}\colon S^1 \to \mathcal{H}$ such that each map $\tilde{\psi}_t$ fixes some interval $\Delta_{j(t)} \in \mathcal{N}^i$.*

**Proof.** We require the following technical definition:

**Definition 5.4.** *We say that a map $f \in \widehat{\mathcal{NS}}$ **transversely immerses** an open subset $\mathcal{I} \subset \gamma$ if $f$ immerses $\mathcal{I}$ and further the subsetset*

$$\mathcal{I}_\cap \subset \mathcal{I} = \{x \in \mathcal{I} \colon \text{there is a } y \in \gamma \text{ such that } f(y) = f(x) \text{ but } x \neq y\}$$

*has an open, dense complement. Let*

$$\mathcal{NS}^i_{\mathrm{tr}} = \{f \in \mathcal{NS}^i \colon \text{ there is a } \Delta_j \in \mathcal{N}^i \text{ transversely immersed by } f\}$$

Then Proposition 5.3 is proved by combining the following two Lemmas. The first, Lemma 5.5, shows that if $h_t$ is homotoped into $\mathcal{NS}^i$, it can be further homotoped so that each map transversely immerses some interval $\Delta_j \in \mathcal{N}^i$. The second, Lemma 5.6, shows that if a map transversely immerses an interval the corresponding homeomorphism in $\mathcal{H}$ must fix it.

**Lemma 5.5.** *The inclusion $i\colon\mathcal{NS}^i_{\mathrm{tr}} \hookrightarrow \mathcal{NS}^i$ is a homotopy equivalence.*

**Proof.** We remind the reader that $\mathcal{U}^j$ denotes the open cover of $\mathcal{NS}^i$ given by:

$$\mathcal{U}^j = \{f \in \mathcal{NS}^i \text{ such that } f|_{\Delta_j} \text{ is an immersion on an interval containing } \Delta_j\}$$

This Lemma is then an easy consequence of the Homotopy Decomposition Lemma (Lemma 3.6), applied to the covering $\mathcal{U}_j$. Let $S \subset J$, $\Delta = \bigcup_{j \in S} \Delta_j$. Then let

$$\mathcal{NS}^\Delta = \bigcap_{j \in S} \mathcal{U}_j$$

be those maps immersing $\Delta$, and let $\mathcal{NS}_{\text{tr}}^\Delta = \mathcal{NS}^\Delta \cap \mathcal{NS}_{\text{tr}}^i$. One can perturb any family of maps in $\mathcal{NS}^\Delta$ to transversely immerse $\Delta$, thus the inclusion $\mathcal{NS}_{\text{tr}}^\Delta \hookrightarrow \mathcal{NS}^\Delta$ is a homotopy equivalence. By the Homotopy Decomposition Lemma, the inclusion $i \colon \mathcal{NS}_{\text{tr}}^i \hookrightarrow \mathcal{NS}^i$ is also a homotopy equivalence.

$\square$

**Lemma 5.6.** *Let $\tilde{h}$ be a loop of maps in $\mathcal{NS}$ homotopic to $h_t$. Let $\tilde{\zeta} \colon S^1 \to \mathcal{H}(H_-, H_+)$ and $\tilde{\psi} \colon S^1 \to \mathcal{H}$ denote loops of homeomorphisms corresponding to $\tilde{h}$. Then if $\tilde{h}_t|_\gamma$ transversely immerses an open interval $\mathcal{I}$, $\tilde{\psi}_t$ fixes each $x \in \mathcal{I}$.*

**Proof.** Let $\mathcal{I}_{\text{em}} = \mathcal{I} \backslash \mathcal{I}_\cap$ (See Definition 5.4). Since $\tilde{h}_t|_\gamma$ transversely immerses $\mathcal{I}$, $\mathcal{I}_{\text{em}} \subset \mathcal{I}$ is dense. We will show that $\psi_t$ fixea $\mathcal{I}_{\text{em}}$; it must then extend to fix all of $\mathcal{I}$ by continuity. Let $x \in \mathcal{I}_{\text{em}}$, $t \in S^1$. Choose a system of half neighborhoods

$$U_-^j \supset U_-^{j+1}...$$

of $x$ within $H_-$, and a corresponding system of half neighborhods

$$U_+^j \supset U_+^{j+1}...$$

within $H_+$, such that

$$\tilde{h}_t(U_+^j) = \tilde{h}_t(U^j)$$

$\tilde{\zeta}_t(U_-^j)$ is a half disc in $U_\zeta^j \subset H_+$ which is a neighborhood of $x$ which adjoins $\mathcal{I}_{\text{em}}$ and further

$$\tilde{h}_t(U_\zeta^j) = \tilde{h}_t(U_-^j) = \tilde{h}_t(U_+^j)$$

By definition $h_t(\mathcal{I}/\mathcal{I}_{\text{em}}) \cap h_t(\mathcal{I}_{\text{em}}) = \varnothing$. So, there is a $j_0 \in \mathbb{N}$ such that, for $j > j_0$, $U_+^j$ is the unique such half disc in $H_+$. Thus $\tilde{\zeta}_t(U_-^j) = U_+^j$, and:

$$\tilde{\psi}_t(U_-^j) = \zeta_0^{-1} \tilde{\zeta}_t(U_-^j)$$

is a half neighborhood of $x$ within $H_-$. Therefore, the sets $\psi_t(U^j)$, for $j > j_0$, give a system of neighborhoods of $x$, and thus $\psi_t$ extends continiously to a map fixing $x$.

$\square$

This completes the proof of Proposition 5.3. For by Lemma 5.5 we can homotop $h$ to a loop $\tilde{h} \colon S^1 \to \mathcal{NS}$ such that each map $\tilde{h}_t$ transversely immerses some interval $\Delta_{j(t)} \in \mathcal{N}^i$. Then is then a corresponding homotopy of $\psi$ to $\tilde{\psi}$. By Lemma 5.6 each map $\tilde{\psi}_t$ then fixes $\Delta_{j(t)}$.

$\square$

## 5.2 Repeated Alexander tricks show $\tilde{\psi}$ is contractible

In this final subsection we complete the proof of Proposition 5.1: *If $[h] \in \pi_1(\widehat{\mathcal{NS}}, \mathcal{NS})$ is trivial then the path $\psi \colon S^1 \to \mathcal{H}$ is contractible.*

We now show that since each $\tilde{\psi}_t$ fixes a segment $\Delta_j \in \mathcal{N}^i$, the loop $\tilde{\psi}$ must be contractible. As $\mathcal{H}$ is a topological group, it is sufficient to show that the corresponding cycle $[\tilde{\psi}]$ is trivial in homology. We will do this through a covering argument. We cannot apply directly the usual statement of Mayer-Veitoris, as the compact open topology on $\mathcal{H}$ is quite weak. However, the loop $\tilde{\psi}$ is continious in a much finer topology -- that induced in some sense from $\mathcal{NS}$. In 5.2.3 we use this finer continuity implicitly to apply the argument behind Mayer-Veitoris and show that $[\tilde{\psi}] = 0$ in $H_1(\mathcal{H})$.

### 5.2.1 Alexander trick shows that homeomorphisms of $D^2$ fixing an arc is contractible

**Lemma 5.7.** *Let* $\Delta \subset \partial D^2$ *denote a closed interval. Then the homeomorphisms* $\mathcal{H}_\Delta \subset \mathcal{H}$ *which fix* $\Delta$ *form a contractible set.*

**Proof.** Let $\rho : (\partial D^k, D^k) \to (\mathrm{id}, \mathcal{H}_\Delta)$ denote a map of pairs. We now construct a deformation retraction $\boldsymbol{\rho} : D^k \times I \to \mathcal{H}_\Delta$ of $\rho$. Consider $D^2$ as a rectangle in $\mathbb{R}^2$:

$$D^2 = \{(x, y) \in \mathbb{R}^2 : 0 < x < 1, 0 < y < 1\}$$

such that:

$$\Delta = \{(x, y) \in \mathbb{R}^2 : 0 < x < 1, y = 1\}$$

Let $D^2_+$ be the double of $D^2$ along $\Delta$:

$$D^2_+ = \{(x, y) \in \mathbb{R}^2 : 0 < x < 1, 0 < y < 2\}$$

We extend each homeomorphism $\rho_d$ to a homeomorphism $\rho_d^+$ of $D^2_+$ by the identity map:

$$\rho_d^+(x) = \begin{cases} \rho_d(x) & \text{if } x \in D^2 \\ x & \text{otherwise} \end{cases}$$

Denote the homeomorphisms of $D^2_+$ by $\mathcal{H}(D^2_+)$. Then we define a path of homeomorphisms $\boldsymbol{\varsigma} : I \to \mathcal{H}(D^2_+)$ by:

$$\varsigma_t(x, y) = \begin{cases} (x, \frac{1}{1-t}y) & \text{if } y < \frac{3}{2}(1-t) \\ (x, \frac{1}{1+3t}y + \frac{6t}{1+3t}) & \text{if } \frac{3}{2}(1-t) \leq y < 2 \end{cases}$$

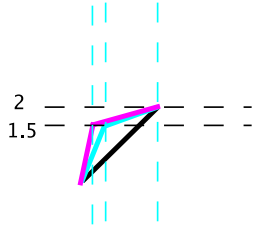

**Figure 5.1.** The graphs of the function $y \to \varsigma_t(\cdot, y)$, as $t$ tends towards 1.

$\varsigma_t$ sends any compact subset of $D^2$ into $D^2_+ \backslash D^2$ for $t$ large enough. So, since $\rho_+|_{D^2_+ \backslash D^2}$ is the identity the deformation retraction $\boldsymbol{\rho} : D^k \times I \to \mathcal{H}_\Delta$ of $\rho$ to the identity given by:

$$\boldsymbol{\rho}_{d,t} = \begin{cases} \varsigma_t^{-1} \circ \rho_d^+ \circ \varsigma_t & 0 \leq t < 1 \\ \mathrm{id} & t = 1 \end{cases}$$

is continous in the topology of convergence on compact sets. $\qquad \square$

### 5.2.2 A contractible covering with contractible intersections

**Definition 5.8.** *Let* $i \in \mathbb{N}$. *We consider the following subsets of* $\mathcal{H}$

1. *Let* $\Delta_j \in \mathcal{N}^i$, *and denote by* $\mathcal{H}_j$ *the homeomorphisms in* $\mathcal{H}$ *which fix* $\Delta_j$.

2. $\mathcal{H}_{j,k} = \mathcal{H}_j \cap \mathcal{H}_k$.

3. $\bar{\mathcal{H}}$ *denotes the homeomorphisms in* $\mathcal{H}$ *which fix the boundary* $\partial D^2$.

**Lemma 5.9.** *Each* $\mathcal{H}_j$ *is contractible and each intersection* $\mathcal{H}_{j,k}$ *is contractible.*

**Proof.** Each $\mathcal{H}_j$ is contractible by Lemma 5.7.

To show that each intersection $\mathcal{H}_{j,k}$ is contractible, fix $y_j \in \Delta_j$, $y_k \in \Delta_k$, and a continious embedding

$$E_\Delta \colon [0,1] \to B^n$$

of the closed interval into $D^2$, such that $E_\Delta(0) = y_j$ and $E_\Delta(1) = y_k$. Let $\mathcal{E}_\Delta$ denote the space of such embedded paths. Then, $\mathcal{E}_\Delta$ is the orbit of $E_\Delta$ under the actions of both $\mathcal{H}_{j,k}$ and $\bar{\mathcal{H}}$. Thus we have the following morphism of fibrations:

$$
\begin{array}{ccc}
\mathcal{H}_{j,k,E} & \to & \bar{\mathcal{H}}_E \\
\downarrow & & \downarrow \\
\mathcal{H}_{j,k} & \to & \bar{\mathcal{H}} \\
\downarrow \phi_{j,k} & & \downarrow \bar\phi \\
\mathcal{E}_\Delta & \to & \mathcal{E}_\Delta
\end{array}
$$

Where $\mathcal{H}_{j,k,E}$ and $\bar{\mathcal{H}}_E$ denote the stabilizers of $E_\Delta$ in their respective groups $\mathcal{H}_{j,k}$ and $\bar{\mathcal{H}}$. We will now show that $\mathcal{H}_{j,k}$ is contractible by showing that every other space in the above diagram is contractible:

**$\bar{\mathcal{H}}$ is contractible:**   We apply the Alexander trick, coning off from the boundary $\partial D^2$.

**$\bar{\mathcal{H}}_E$ is contractible:**   $D^2 \backslash E_\Delta$ consists of two components each homeomorphic to $D^2$. Restricting each homeomorphisms $f \in \bar{\mathcal{H}}_E$ to this pair of discs yields a homeomorphism:

$$\bar{\mathcal{H}}_E \to \bar{\mathcal{H}} \times \bar{\mathcal{H}}$$

$$f \to f|_{D^2 \backslash E_\Delta}$$

Thus $\bar{\mathcal{H}}_E$ is also contractible.

**$\mathcal{E}_\Delta$ is contractible:**   The orbit map $\bar\phi \colon \bar{\mathcal{H}} \to \mathcal{E}_\Delta$ has contractible fiber $\bar{\mathcal{H}}_E$ and thus it is a homotopy equivalence. So, $\mathcal{E}_\Delta$ is contractible as $\bar{\mathcal{H}}$ is contractible.

**$\mathcal{H}_{j,k,E}$ is contractible:**   Again, restricting each homeomorphisms $f \in \mathcal{H}_{j,k,E}$ to the pair of discs $D^2 \backslash E_\Delta$ yields a homeomorphism:

$$\bar{\mathcal{H}}_E \to \mathcal{H}_\Delta \times \mathcal{H}_\Delta$$

$$f \to f|_{D^2 \backslash E_\Delta}$$

and thus $\mathcal{H}_{j,k,E}$ is contractible since each factor $\mathcal{H}_\Delta$ is contractible by Lemma 5.7.

Finally, we see that:

**$\mathcal{H}_{j,k}$ is contractible:**   Since its fiber $\mathcal{H}_{j,k,E}$ is contractible, the orbit map $\phi_{j,k} \colon \mathcal{H}_{j,k} \to \mathcal{E}_\Delta$ is a homotopy equivalence. Therefore as $\mathcal{E}_\Delta$ is contractible, so is $\mathcal{H}_{j,k}$. This completes the proof of Lemma 5.9.   □

### 5.2.3   Mayer − Vietoris type proof that $[\psi]$ is trivial in $H_1(\mathcal{H})$

We now apply the argument behind Mayer-Vietoris to show that $[\psi]$ is trivial in $H_1(\mathcal{H})$.

**Definition 5.10.**  *For any $i \in \mathbb{N}$, let $\mathcal{V}_j \subset \mathcal{NS}^i$ denote those maps which transversely immerse $\Delta_j \in \mathcal{N}^i$ Then*

$$\mathcal{NS}^i_{\mathrm{tr}} = \bigcup_{j \in J} \mathcal{V}_j$$

*and the sets $\mathcal{V}_j$ then give an open cover of $\mathcal{NS}^i_{\mathrm{tr}}$.*

Consider a set of points $x_1, ..., x_n \in S^1$ such that for each interval $[x_i, x_{i+1}]$ there is some set $\mathcal{V}_j$ with $\tilde{h}_{[x_i, x_{i+1}]} \subset \mathcal{V}_j$. [5.1] The corresponding homeomorphisms $\tilde\psi_{[x_i, x_{i+1}]}$ then lie in $\mathcal{H}_j$ by Proposition 5.6.

---

5.1. Here, and in the remainder of this argument, we consider the indices of the $x_i$ as lying in $\mathbb{Z}_n$, and thus consider the interval $[x_n, x_1]$ as an interval of type $[x_i, x_{i+1}]$.

Suppose $\tilde{h}_{[x_i, x_{i+1}]} \subset \mathcal{V}_j$ and $\tilde{h}_{[x_{i-1}, x_i]} \subset \mathcal{V}_k$. Then $\psi_{x_i} \in \mathcal{H}_{j,k} = \mathcal{H}_j \cap \mathcal{H}_k$. So, as $\mathcal{H}_{j,k}$ is connected by Lemma 5.9, we can find a path $p_i \colon I \to \mathcal{H}_{j,k}$ connecting $\psi_{x_i}$ to the identity. Thus, we can write the cycle $[\psi]$ as a sum of cycles:

$$[\psi] = \sum_{i \in \mathbb{Z}_n} [p_i] + [\psi_{[x_i, x_{i+1}]}] - [p_{i+1}]$$

where each summand:

$$[p_i] + [\psi_{[x_i, x_{i+1}]}] - [p_{i+1}]$$

gives a cycle in one of the subsets $\mathcal{H}_j$. Since $H_1(\mathcal{H}_j) = 0$ (again by Lemma 5.9), $[\psi]$ is a sum of null homologous cycles and is thus null homologous.

**Remark 5.11. (Extending this example to higher dimensions)** Most of the proof of Theorem 2.6 goes through without change to dimensions $n > 2$. The construction of model spaces $\mathcal{NS}$ and $\widehat{\mathcal{NS}}$ is the same in higher dimensions. One can still define an obstruction cycle $h \in \pi_1(\mathcal{NS}, \widehat{\mathcal{NS}})$, rotating $H_+$ against $H_-$ in the same way as in dimension 2, although to do so one must now require the $nq$ jets of the forms $\mu_\kappa$ to vanish along the hypersurface $\gamma$. Our argument still shows that the sets $\mathcal{H}_j$ are contractible. *However the author does not know how to generalize the proof that their intersections $\mathcal{H}_{j,k}$ are connected (our argument requires only connectivity) to higher dimensions.*

# 6 Appendix

## 6.1 Overlapping families of maps

In this subsection we prove Lemma 3.7.

**Lemma. (3.7)** *Let $\gamma_{em} \subset \gamma$ be a closed interval, let $U \subset D^k$ be an open set, and let $\rho \colon U \to \mathcal{NS}_{nc}$ is a family of maps which each embedd $\gamma_{em}$. Then there is a homotopy $\boldsymbol{\rho} \colon U_\Delta \times I \to \mathcal{NS}_{nc}$ such that:*

1. $\boldsymbol{\rho}_{d,0} = \rho_d$

2. $\boldsymbol{\rho}_{d,1}$ *overlaps, and thus* $\boldsymbol{\rho}_{d,1} \in \mathcal{NS}$.

3. *If $\boldsymbol{\rho}_{d,0}$ overlaps, then $\boldsymbol{\rho}_{d,t}$ overlaps for all $t \in [0,1]$.*

4. *If $\boldsymbol{\rho}_{d,0}$ immerses an interval $\Delta \subset \gamma$ then $\boldsymbol{\rho}_{d,t}$ immerses $\Delta$ for all $t \in [0,1]$.*

**Proof.** We begin by choosing a family of embedded bands $\boldsymbol{B}$ such that $\boldsymbol{B}_d$ is a thickenings of an arc which meet $\rho_d(\gamma_{em})$ transversely. Then we push each map along these bands so that they eventually overlap. Finally we adjust this pushing deformation so that we still induce the proper form near $\gamma$ by reparemtarizing in the domain.

Let $V_{em}$ be a neighborhood of $\gamma_{em}$ such that both $V_{em} \cap H_+$ and $V_{em} \cap H_-$ are embedded by each $\rho_d$.

**Proposition 6.1.** *Denote by $\mathcal{E}mb(D^2, S^2)$ the space of embeddings of $D^2$ into $S^2$. Then there is a continuous family of embeddings $B \colon D^k \to \mathcal{E}mb(I_\delta, S^2)$, along with a pair of disjoint discs $I^r$ and $I^l$ inside $D^2$ such that:*

1. $B_d(I^l)$ *and* $B_d(I^r)$ *are both contained in $\rho_d(V_{em})$.*

2. $B_d^{-1}(\rho_d(\gamma_{em}))$ *is a pair of arcs such that $D^2 \backslash B_d^{-1}(\rho_d(\gamma_{em}))$ consists of three components: $D_l$, $D_m$, $D_r$. Further $D_l, D_r \subset \rho_d(V_{em})$, and $D_m$ is contained in its complement.*

3. *The closure of* $\operatorname{im}(B_d) \cup \operatorname{im}(\rho_d) \neq S^2$

**Proof.** Fix a metric $g$ on the range $M$. Then the exponential map of $g$ induces a family of embeddings of

$$\mathrm{Exp}\colon D^k \to \mathcal{E}mb((-\varepsilon, \varepsilon) \times \gamma_{\mathrm{em}}, M)$$

such that

$$\mathrm{Exp}_d(0, x) = \rho_d(x) \text{ for each } x \in \gamma_{\mathrm{em}}$$

Now choose one embedding

$$f\colon [-1, 1] \times [-1, 1] \to [-\varepsilon, \varepsilon] \times \gamma_{\mathrm{em}}$$

such that

$$f^{-1}([-\varepsilon, 0] \times \gamma_{\mathrm{em}}) = [-1, -\tfrac{1}{2}] \times [-1, 1] \cup [\tfrac{1}{2}, 1] \times [-1, 1]$$

Then for sufficiently small $\delta > 0$:

$$\mathrm{Exp}_d \circ f([-1, 1] \times [-\delta, \delta]) \cup \mathrm{im}\,(\rho_d) \neq S^2$$

for all $d \in D^k$. Thus

$$B_d = \mathrm{Exp}_d \circ f|_{[-1,1] \times [-\delta, \delta]}$$

gives the required family of embeddings.

$\square$

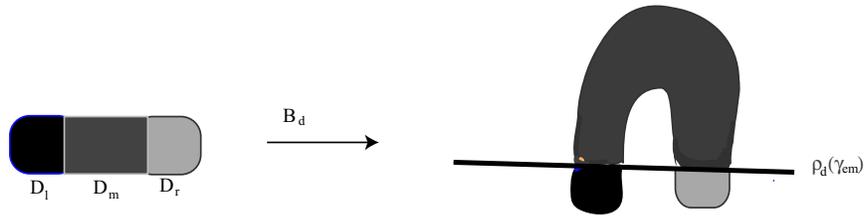

**Figure 6.1.** The emmbeded bands $B_d$. The shaded regions represent the sets $D_l$, $D_m$ and $D_r$ on the left, and their images under $B_d$ on the right.

Let $\varsigma\colon (0, 1) \to \mathrm{Diff}(D^2)$ be a continuous family of diffeomorphisms of $D^2$ such that:

1. $\varsigma_0$ is the identity map.

2. $\varsigma_t$ fixes the boundary of $D^2$.

3. For $t \geq \tfrac{1}{2}$,

$$\varsigma_t(D_r) \cap D_l \neq \varnothing \tag{6.1}$$

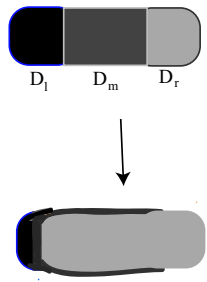

**Figure 6.2.** The diffeomorphism $\varsigma_1$.

Now let $\mathcal{NS}'_{\mathrm{nc}}$ denote the degree 0, $C^q$, non-surjective, maps which are local diffeomorphisms away from $\gamma$, orientation preserving in $H_+$ and reversing in $H_-$. Then $\mathcal{NS}'_{\mathrm{nc}} \supset \mathcal{NS}_{\mathrm{nc}}$, and *the two differ only in that a map $f \in \mathcal{NS}'_{\mathrm{nc}}$ need not satisfy $f^*\sigma = \mu_1$ near $\gamma$.* Similarly denote by $\mathcal{NS}' \subset \mathcal{NS}'_{\mathrm{nc}}$ those maps $f \in \mathcal{NS}'_{\mathrm{nc}}$ which overlap.

Let

$$U_d^l = \rho_d^{-1}(B_d(D_l)) \cap V_{\mathrm{em}}$$
$$U_d^r = \rho_d^{-1}(B_d(D_r)) \cap V_{\mathrm{em}}$$

Then in order to define $\boldsymbol{\rho} \colon U_\Delta \times I \to \mathcal{NS}_{\mathrm{nc}}$ overlapping $\rho$, we first define

$$\boldsymbol{\rho}' \colon U_\Delta \times I \to \mathcal{NS}'_{\mathrm{nc}}$$

by:

$$\boldsymbol{\rho}'_{d,t}(x) = \begin{cases} x & \text{if } x \in S \setminus U_d^r \\ B_d \circ \varsigma_t \circ B_d^{-1} \circ \rho_d(x) & x \in U_d^r \end{cases}$$

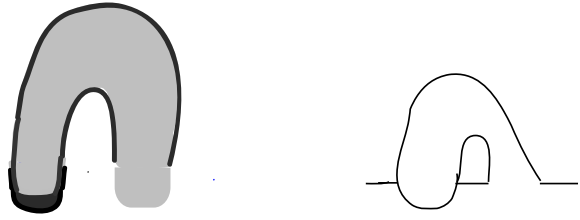

**Figure 6.3.** The image of $B_d \circ \varsigma_1$ (left), and $\boldsymbol{\rho}'_{d,1}$ (right).

We now verify that, while $\boldsymbol{\rho}'_{d,t}{}^*\sigma \neq \mu_1$ near $\gamma$, it satisfies the other requirements of our Lemma.

1. $\boldsymbol{\rho}'_{d,t} \in \mathcal{NS}'_{\mathrm{nc}}$ for all $d$, $t$: $\boldsymbol{\rho}'_{d,t}$ is given by composing $\rho_d$ with local diffeomorphisms. Therefore $\boldsymbol{\rho}'_{d,t}$ has the proper $q - $ jets away from $\gamma$, as $\rho_d$ did. Finally the deformation $\boldsymbol{\rho}'_{d,t}$ proceeds through non-surjective maps. For the $\mathrm{im}(\boldsymbol{\rho}'_{d,t}) \subset \mathrm{im}(B_d) \cup \mathrm{im}(\rho_d)$, and the bands $B_d$ are constructed to satisfy the condition

$$\mathrm{im}(B_d) \cup \mathrm{im}(\rho_d) \neq S^2$$

2. $\boldsymbol{\rho}'_{d,1} \in \mathcal{NS}'$: We need to show $\boldsymbol{\rho}'_{d,1}$ overlaps.

   I claim that

   $$\boldsymbol{\rho}'_{d,1}(U_d^r) \cap \boldsymbol{\rho}'_{d,1}(U_d^l) \neq \varnothing$$

   and thus $\boldsymbol{\rho}'_{d,1}$ overlaps. For $\boldsymbol{\rho}'_{d,1}(U_d^l) = \rho_d(U_d^l) = B_d(D_l)$, and $\boldsymbol{\rho}'_{d,1}(U_d^r) = B_d(\varsigma_1 \circ D_r)$. Then

   $$\boldsymbol{\rho}'_{d,1}(U_d^r) \cap \boldsymbol{\rho}'_{d,1}(U_d^l) = B_d(D_l) \cap B_d(\varsigma_1 \circ D_r)$$
   $$= B_d(D_l \cap \varsigma_1 \circ D_r)$$

   which is non-empty as $D_l \cap \varsigma_1 \circ D_r$ is nonempty. (equation 6.1).

3. If $\rho_d$ overlaps, then $\boldsymbol{\rho}'_{d,t}$ overlaps for all $t$: $\boldsymbol{\rho}'_{d,t} = \rho_d$ outside of $U_d^r$, and $\boldsymbol{\rho}'_{d,t}(U_d^r) \supset \rho_d(U_d^r)$. Thus $\boldsymbol{\rho}'_{d,t}|_{H_+}$ is non-injective if $\rho_d|_{H_+}$ was.

4. If $\rho_d$ immerses $\Delta$ then so does $\boldsymbol{\rho}'_{d,t}$, as the latter is given by composing the former with local diffeomorphisms.

$\boldsymbol{\rho}'$ is very nearly the "overlapping" deformation we require. However we must make a final adjustment, for in general $\boldsymbol{\rho}'_{d,t}{}^*\sigma \neq \mu_1$ near $\gamma$, and thus $\boldsymbol{\rho}'_{d,t}$ doesn't remain in $\mathcal{NS}_{\mathrm{nc}}$, but lies only in $\mathcal{NS}'_{\mathrm{nc}}$. We correct this by precomposing with a diffeomorphism of the domain.

**Modifying our deformation $\boldsymbol{\rho}'$ to preserve $\mu_1$ near $\gamma$**

Let:

$$K_\rho(d,x,t) = \begin{cases} 1 & \text{if } x \in S \setminus U_d^r \\ \dfrac{(B_d \circ \varsigma_{\phi(d) \cdot t} \circ B_d^{-1})^* \sigma(\rho_d(x))}{\sigma(\rho_d(x))} & x \in U_d^r \end{cases}$$

$K_\rho(x,t)$ is continuous in both variables.

Let $(\theta, r)$ be sphereical coordinates on the domain $M = S^2$, so that $\theta \in [0, 2\pi)$, $r \in [-1, 1]$ and $\gamma = \{x: r(x) = 0\}$. Let $g_t$ be a diffeomorphism of the domain such that near $\gamma$:

$$\boldsymbol{g}_{d,t}(\theta, h) = (\theta, \frac{1}{K_\rho(d,x,t)} r)$$

Then we define the homotopy $\boldsymbol{\rho}: D^k \times I \to \mathcal{N}\mathcal{S}_{\mathrm{nc}}$ by:

$$\boldsymbol{\rho}_{d,t} = \boldsymbol{\rho}'_{d,t} \circ \boldsymbol{g}_{d,t}$$

then this deformation is trivial away from each $U_r(d)$ and inside these we have:

$$\begin{aligned}
\boldsymbol{\rho}^*_{d,t}(\sigma) &= \boldsymbol{g}^*_{d,t}(\boldsymbol{\rho}'_{d,t}{}^*\sigma) \\
&= \boldsymbol{g}^*_{d,t}\boldsymbol{\rho}^*_d\big((B_d \circ \varsigma_{\phi(d) \cdot t} \circ B^{-1}{}_d)^*\sigma\big) \\
&= \boldsymbol{g}^*_{d,t}(K_\rho(d,x,t)\sigma) \\
&= \mu_1
\end{aligned}$$

and so each map $\boldsymbol{\rho}_{d,t}$ induces the proper form near $\gamma$, and thus $\boldsymbol{\rho}_{d,t} \in \mathcal{N}\mathcal{S}_{\mathrm{nc}}$ for all $d$, $t$. The remaining conditions of our claim continue to hold, as $\boldsymbol{\rho}_{d,t}$ differs from $\boldsymbol{\rho}'_{d,t}$ only by a reparameterization of the domain.

This completes the proof of Lemma 3.7.                                                $\square$

## 6.2   Proof of Gray's Lemma

In this subsection we provide the necessary induction to prove Gray's Lemma (Lemma 3.6).

**Lemma. (Homotopy Decomposition Lemma [Gra75] - Proposition 16.24)**   *Let* $f: X \to Y$ *be a continuous map. Let* $U_Y$ *be a finite covering of* $Y$ *by open sets* $U_Y^j$, *and denote* $f^{-1}(U_Y^j)$ *by* $U_X^j$. *Suppose that for each* $J \subset I$ *the the restriction*

$$f: \bigcap_{j \in J} U_X^j \to \bigcap_{j \in J} U_Y^j$$

*is a homotopy equivalence then* $f$ *is a homotopy equivalence*

**Proof.**   Proposition 16.24 in [Gra75] covers the case of a covering by 2 sets. The general case follows by an induction. Let $U_Y'$ be the refinement of $U_Y$ given by all of its multi-intersections, and let $U_X'$ be the analagous refinement of $U_X$. Then the mulintersections of the refinements $U_X'$ and $U_Y'$ are the same as those of the original covers $U_X$ and $U_Y$, and so we see that $f$ restricted to each multi-intersection of $U_X'$ is a homopty equivalence.

Suppose that we have shown that $f$ is a homotopy equivalence  when restricted oto any set $U'$ given by a union of $l$ members of $U_X'$.  Then I claim that

$$f|_{U' \cup U_X^{jk}}$$

is a homotopy equivalence for any $U_X^{jk}$ in the cover $U_X'$. For $U' \cap U_X^{jk}$ can also be written as the union of $l$ members of $U_X'$. Thus, by our induction hypothesis,  $f|_{U' \cap U_X^{jk}}$ is a homotopy equivalence.  So by Gray's Proposition 16.24, applied to the sets $U'$ and $U_X^{jk}$ we see that $f|_{U' \cup U_X^{jk}}$ is also a homotopy equivalence.   $\square$

## 6.3   Constructing scaling functions

In this subsection we provide the proofs of Lemmas  3.10 and 3.11, the construction of the functions which scale our forms in the image, and co-image of the maps $\rho_d$.

**Lemma. (3.10)** *For each $0 < \delta$, there is a continious function $\phi_{\mathrm{im}}^{\delta}\colon D^k \times S^2 \to [0, \infty)$ such that each restriction $\phi_{\mathrm{im}}^{\delta}|_{d \times S^2}$ is $C^{\infty}$ and $V_d^{\delta} = (\phi_{\mathrm{im}}^{\delta})^{-1}(0) \cap (d \times S^2)$ satisfies:*

1. *$\rho_d^{-1}(V_d^{\delta})$ is a closed neighborhood of $\gamma$ such that:*

$$\int_{\rho_d^{-1}(V_d^{\delta})\ \cap H_+} \rho_d^* \sigma < \delta$$

2. *$V_d^{\delta} \cup \mathrm{im}(\rho_d) = S^2$.*

**Proof.** Choose an auxilliary metric on the range $N$. Define a continuous function $\lambda\colon D^k \to [0, \infty)$ such that:

1. $\lambda_d = 0$ if and only if $\rho_d$ is surjective.

2. If $U_d^{\lambda}$ is the closed $\lambda_d$ neighborhood of $S^2 \backslash \mathrm{im}(\rho_d)$, then

$$\int_{\rho_d^{-1}(U_d^{\lambda}) \cap \mathrm{im}(\rho_d)} \sigma < \delta$$

Let $\phi_{\mathrm{cont}}\colon D^k \times S^2 \to [0, \infty)$ denote a continious function such that:

$$\phi_{\mathrm{cont}}^{-1}(0) \cap d \times S^2 = U_d^{\varepsilon}$$

We now smooth $\phi_{\mathrm{cont}}$ in the $S^2$ direction, by convoluting with a continious family of bump functions $\varphi$, parameterized by $D^k$, and such that $\varphi_d$ is supported in a ball of radius $\frac{\lambda_d}{4}$, and also of height $\frac{\lambda_d}{4}$. We thus gain our function $\phi_{\mathrm{im}}^{\delta}\colon D^k \times S^2 \to [0, \infty)$, which is smooth along each $d \times S^2$.

I claim $\phi_{\mathrm{im}}^{\delta}$ satisfies the conditions of the Lemma: $V_d^{\delta}$ consists of the points $x \in S^2$ such that an $\frac{\lambda_d}{4}$ ball is contained in $U_d^{\lambda}$. The triangle inequality then ensures that $V_d^{\delta}$ contains an $\frac{\lambda_d}{4}$ neighborhood of $S^2 \backslash \mathrm{im}(\rho_d)$. In particular $\rho_d^{-1}(V_d^{\delta})$ is a closed neighborhood of $\gamma$, and $V_d^{\delta} \cup \mathrm{im}(\rho_d) = S^2$. Moreover, since $V_d^{\delta} \subset U_d^{\lambda}$ we have that:

$$\int_{\rho_d^{-1}(V_d^{\delta})\ \cap H_+} \rho_d^* \sigma(x) \leq \int_{\rho_d^{-1}(U_d^{\lambda}) \cap \mathrm{im}(\rho_d)} \sigma < \delta \qquad\qquad \square$$

We now prove Lemma 3.11, the construction of the second scaling function.

**Lemma. (3.11)** *Let $\varepsilon\colon D^k \to [0, \infty)$ be a continious function such that $\varepsilon_d = 0$ if and only if $\rho_d$ is surjective. Then there is a continious function $\phi_{\mathrm{ms}}^{\varepsilon}\colon D^k \times S^2 \to [0, -\infty)$ such that:*

1. *Each restriction $\phi_{\mathrm{ms}}^{\varepsilon}|_{d \times S^2}$ is $C^{\infty}$.*

2. *$V_d^{\varepsilon} = \phi_{\mathrm{ms}}^{\varepsilon}{}^{-1}(0) \cap (d \times S^2)$ is a closed neighborhood of $\mathrm{im}(\rho_d)$ within $s \times S^2$ satisfying:*

$$\int_{V_d^{\varepsilon}} \nu_d' - \int_{\mathrm{im}(\rho_d)} \nu_d' \leq \varepsilon_d$$

**Proof.** Define a continuous function $\lambda\colon D^k \to [0, \infty)$ such that:

1. $\lambda_d = 0$ if and only if $\rho_d$ is surjective.

2. If $U_d^{\lambda}$ is the closed $\lambda_d$ neighborhood of $\mathrm{im}(\rho_d)$, then

$$\int_{U_d^{\lambda}} \nu_d' - \int_{\mathrm{im}(\rho_d)} \nu_d' \leq \varepsilon_d$$

Let $\phi_{\mathrm{cont}}\colon D^k \times S^2 \to [0, \infty)$ denote a continious function such that:

$$\phi_{\mathrm{cont}}^{-1}(0) \cap d \times S^2 = U_d^{\lambda}$$

We now smooth $\phi_{\text{cont}}$ in the $S^2$ direction, by convoluting with a family of bump functions, each supported in a ball of radius $\frac{\lambda_d}{4}$, to gain our function $\phi_{\text{ms}}^\varepsilon \colon D^k \times S^2 \to [0, \infty)$.

I claim $\phi_{\text{ms}}^\varepsilon$ satisfies the conditions of the Lemma: $V_d^\varepsilon$ consists of the points $x \in S^2$ such that an $\frac{\lambda_d}{4}$ ball is contained in $U_d^\lambda$. The triangle inequality then ensures that $V_d^\varepsilon$ contains an $\frac{\lambda_d}{4}$ neighborhood of $\operatorname{im}(\rho_d)$. Moreover, since $V_d^\varepsilon \subset U_d^\lambda$ we have that:

$$\int_{\rho_d^{-1}(V_d^\varepsilon) \,\cap H_+} \rho_d^* \sigma(x) \leq \int_{\rho_d^{-1}(U_d^\lambda) \cap \operatorname{im}(\rho_d)} \sigma \leq \varepsilon_d \qquad\qquad \square$$